\documentclass[11pt]{article}
\usepackage[dvips]{graphicx}
\usepackage{textcomp}
\usepackage{amsbsy}
\usepackage{latexsym}
\usepackage[mathscr]{eucal}
\usepackage{amsfonts,amsmath,amsthm}
\usepackage{amssymb}
\usepackage[usenames]{xcolor}

\usepackage{fullpage}
\def\anew{\color{black} }

\def\cnew{\color{black} }
\def\RIP{\mathop{\rm RIP}}
\begin{document}
\bibliographystyle{plain}
\title
{Optimal Stable Nonlinear Approximation}
\author{ 
Albert Cohen, Ronald DeVore,   Guergana Petrova, and Przemyslaw Wojtaszczyk
\thanks{%
   This research was supported by the  NSF Grant  DMS 1817603 (RD-GP) and the ONR Contract N00014-17-1-2908 (RD).   P. W. was supported by National Science Centre, 
   Polish grant UMO-2016/21/B/ST1/00241. A portion of this research was
   completed when the first three authors were visiting the Isaac Newton Institute. }}
\hbadness=10000
\vbadness=10000
\newtheorem{lemma}{Lemma}[section]
\newtheorem{prop}[lemma]{Proposition}
\newtheorem{cor}[lemma]{Corollary}
\newtheorem{theorem}[lemma]{Theorem}
\newtheorem{remark}[lemma]{Remark}
\newtheorem{example}[lemma]{Example}
\newtheorem{definition}[lemma]{Definition}
\newtheorem{proper}[lemma]{Properties}
\newtheorem{assumption}[lemma]{Assumption}
%
\newenvironment{disarray}{\everymath{\displaystyle\everymath{}}\array}{\endarray}

\def\RR{\rm \hbox{I\kern-.2em\hbox{R}}}
\def\NN{\rm \hbox{I\kern-.2em\hbox{N}}}
\def\ZZ{\rm {{\rm Z}\kern-.28em{\rm Z}}}
\def\CC{\rm \hbox{C\kern -.5em {\raise .32ex \hbox{$\scriptscriptstyle
|$}}\kern
-.22em{\raise .6ex \hbox{$\scriptscriptstyle |$}}\kern .4em}}
\def\vp{\varphi}
\def\<{\langle}
\def\>{\rangle}
\def\t{\tilde}
\def\i{\infty}
\def\e{\varepsilon}
\def\sm{\setminus}
\def\nl{\newline}
\def\o{\overline}
\def\wt{\widetilde}
\def\wh{\widehat}
\def\cT{{\cal T}}
\def\cA{{\cal A}}
\def\cI{{\cal I}}
\def\cV{{\cal V}}
\def\cB{{\cal B}}
\def\cF{{\cal F}}
\def\cY{{\cal Y}}

\def\cD{{\cal D}}
\def\cP{{\cal P}}
\def\cJ{{\cal J}}
\def\cM{{\cal M}}
\def\cO{{\cal O}}
\def\Chi{\raise .3ex
\hbox{\large $\chi$}} \def\vp{\varphi}
\def\lsima{\hbox{\kern -.6em\raisebox{-1ex}{$~\stackrel{\textstyle<}{\sim}~$}}\kern -.4em}
\def\lsim{\hbox{\kern -.2em\raisebox{-1ex}{$~\stackrel{\textstyle<}{\sim}~$}}\kern -.2em}
\def\[{\Bigl [}
\def\]{\Bigr ]}
\def\({\Bigl (}
\def\){\Bigr )}
\def\[{\Bigl [}
\def\]{\Bigr ]}
\def\({\Bigl (}
\def\){\Bigr )}
\def\L{\pounds}
\def\pr{{\rm Prob}}
\newcommand{\cs}[1]{{\color{magenta}{#1}}}
\def\ds{\displaystyle}
\def\ev#1{\vec{#1}}     
\newcommand{\lt}{\ell_{2}(\nabla)}
\def\Supp#1{{\rm supp\,}{#1}}
\def\R{\mathbb{R}}
\def\E{\mathbb{E}}
\def\nl{\newline}
\def\T{{\relax\ifmmode I\!\!\hspace{-1pt}T\else$I\!\!\hspace{-1pt}T$\fi}}
\def\N{\mathbb{N}}
\def\Z{\mathbb{Z}}
\def\N{\mathbb{N}}
\def\Zd{\Z^d}
\def\Q{\mathbb{Q}}
\def\C{\mathbb{C}}
\def\Rd{\R^d}
\def\gsim{\mathrel{\raisebox{-4pt}{$\stackrel{\textstyle>}{\sim}$}}}
\def\sime{\raisebox{0ex}{$~\stackrel{\textstyle\sim}{=}~$}}
\def\lsim{\raisebox{-1ex}{$~\stackrel{\textstyle<}{\sim}~$}}
\def\div{\mbox{ div }}
\def\M{M}  \def\NN{N}                  
\def\Le{{\ell_1}}            
\def\Lz{{\ell_2}}
\def\Let{{\tilde\ell_1}}     
\def\Lzt{{\tilde\ell_2}}
\def\Ltw{\ell_\tau^w(\nabla)}
\def\t#1{\tilde{#1}}
\def\la{\lambda}
\def\La{\Lambda}
\def\ga{\gamma}
\def\BV{{\rm BV}}
\def\Ga{\eta}
\def\al{\alpha}
\def\cZ{{\cal Z}}
\def\cA{{\cal A}}
\def\cU{{\cal U}}
\def\argmin{\mathop{\rm argmin}}
\def\argmax{\mathop{\rm argmax}}
\def\prob{\mathop{\rm prob}}

\def\cO{{\cal O}}
\def\cA{{\cal A}}
\def\cC{{\cal C}}
\def\cS{{\cal F}}
\def\bu{{\bf u}}
\def\bz{{\bf z}}
\def\bZ{{\bf Z}}
\def\bI{{\bf I}}
\def\cE{{\cal E}}
\def\cD{{\cal D}}
\def\cG{{\cal G}}
\def\cI{{\cal I}}
\def\cJ{{\cal J}}
\def\cM{{\cal M}}
\def\cN{{\cal N}}
\def\cT{{\cal T}}
\def\cU{{\cal U}}
\def\cV{{\cal V}}
\def\cW{{\cal W}}
\def\cL{{\cal L}}
\def\cB{{\cal B}}
\def\cG{{\cal G}}
\def\cK{{\cal K}}
\def\cX{{\cal X}}
\def\cS{{\cal S}}
\def\cP{{\cal P}}
\def\cQ{{\cal Q}}
\def\cR{{\cal R}}
\def\cU{{\cal U}}
\def\bL{{\bf L}}
\def\bl{{\bf l}}
\def\bK{{\bf K}}
\def\bC{{\bf C}}
\def\X{X\in\{L,R\}}
\def\ph{{\varphi}}
\def\D{{\Delta}}
\def\H{{\cal H}}
\def\bM{{\bf M}}
\def\bx{{\bf x}}
\def\bj{{\bf j}}
\def\bG{{\bf G}}
\def\bP{{\bf P}}
\def\bW{{\bf W}}
\def\bT{{\bf T}}
\def\bV{{\bf V}}
\def\bv{{\bf v}}
\def\bt{{\bf t}}
\def\bz{{\bf z}}
\def\bw{{\bf w}}
\def \span{{\rm span}}
\def \meas {{\rm meas}}
\def\rhom{{\rho^m}}
\def\diff{\hbox{\tiny $\Delta$}}
\def\EE{{\rm Exp}}
\def\lll{\langle}
\def\argmin{\mathop{\rm argmin}}
\def\codim{\mathop{\rm codim}}
\def\rank{\mathop{\rm rank}}

\def\argmax{\mathop{\rm argmax}}
\def\dJ{\nabla}
\newcommand{\ba}{{\bf a}}
\newcommand{\bb}{{\bf b}}
\newcommand{\bc}{{\bf c}}
\newcommand{\bd}{{\bf d}}
\newcommand{\bs}{{\bf s}}
\newcommand{\bff}{{\bf f}}
\newcommand{\bp}{{\bf p}}
\newcommand{\bg}{{\bf g}}
\newcommand{\by}{{\bf y}}
\newcommand{\br}{{\bf r}}
\newcommand{\be}{\begin{equation}}
\newcommand{\ee}{\end{equation}}
\newcommand{\bea}{$$ \begin{array}{lll}}
\newcommand{\eea}{\end{array} $$}
\def \Vol{\mathop{\rm  Vol}}
\def \mes{\mathop{\rm mes}}
\def \Prob{\mathop{\rm  Prob}}
\def \exp{\mathop{\rm    exp}}
\def \sign{\mathop{\rm   sign}}
\def \sp{\mathop{\rm   span}}
\def \rad{\mathop{\rm   rad}}
\def \vphi{{\varphi}}
\def \csp{\overline \mathop{\rm   span}}
%
%
\newcommand{\beqn}{\begin{equation}}
\newcommand{\eeqn}{\end{equation}}
\def\beginproof{\noindent{\bf Proof:}~ }
\def\endproof{\hfill\rule{1.5mm}{1.5mm}\\[2mm]}

\newenvironment{Proof}{\noindent{\bf Proof:}\quad}{\endproof}

\renewcommand{\theequation}{\thesection.\arabic{equation}}
\renewcommand{\thefigure}{\thesection.\arabic{figure}}

\makeatletter
\@addtoreset{equation}{section}
\makeatother

\newcommand\abs[1]{\left|#1\right|}
\newcommand\clos{\mathop{\rm clos}\nolimits}
\newcommand\trunc{\mathop{\rm trunc}\nolimits}
\renewcommand\d{d}
\newcommand\dd{d}
\newcommand\diag{\mathop{\rm diag}}
\newcommand\dist{\mathop{\rm dist}}
\newcommand\diam{\mathop{\rm diam}}
\newcommand\cond{\mathop{\rm cond}\nolimits}
\newcommand\eref[1]{{\rm (\ref{#1})}}
\newcommand{\iref}[1]{{\rm (\ref{#1})}}
\newcommand\Hnorm[1]{\norm{#1}_{H^s([0,1])}}
\def\int{\intop\limits}
\renewcommand\labelenumi{(\roman{enumi})}
\newcommand\lnorm[1]{\norm{#1}_{\ell^2(\Z)}}
\newcommand\Lnorm[1]{\norm{#1}_{L_2([0,1])}}
\newcommand\LR{{L_2(\R)}}
\newcommand\LRnorm[1]{\norm{#1}_\LR}
\newcommand\Matrix[2]{\hphantom{#1}_#2#1}
\newcommand\norm[1]{\left\|#1\right\|}
\newcommand\ogauss[1]{\left\lceil#1\right\rceil}
\newcommand{\QED}{\hfill
\raisebox{-2pt}{\rule{5.6pt}{8pt}\rule{4pt}{0pt}}%
  \smallskip\par}
\newcommand\Rscalar[1]{\scalar{#1}_\R}
\newcommand\scalar[1]{\left(#1\right)}
\newcommand\Scalar[1]{\scalar{#1}_{[0,1]}}
\newcommand\Span{\mathop{\rm span}}
\newcommand\supp{\mathop{\rm supp}}
\newcommand\ugauss[1]{\left\lfloor#1\right\rfloor}
\newcommand\with{\, : \,}
\newcommand\Null{{\bf 0}}
\newcommand\bA{{\bf A}}
\newcommand\bB{{\bf B}}
\newcommand\bR{{\bf R}}
\newcommand\bD{{\bf D}}
\newcommand\bE{{\bf E}}
\newcommand\bF{{\bf F}}
\newcommand\bH{{\bf H}}
\newcommand\bU{{\bf U}}
\newcommand \A {{\bb A}}
\newcommand\cH{{\cal H}}
\newcommand\sinc{{\rm sinc}}
\def\enorm#1{| \! | \! | #1 | \! | \! |}

\newcommand{\am}{a_{\min}}
\newcommand{\aM}{a_{\max}}

\newcommand{\dm}{\frac{d-1}{d}}

\let\bm\bf
\newcommand{\bbeta}{{\mbox{\boldmath$\beta$}}}
\newcommand{\bal}{{\mbox{\boldmath$\alpha$}}}
\newcommand{\bbi}{{\bm i}}

\def\nnew{\color{Red}}
\def\mnew{\color{Blue}}
\def\wnew{\color{magenta}}

\newcommand{\dI}{\Delta}
\newcommand\aconv{\mathop{\rm absconv}}

\maketitle
 %
  \begin{abstract}   While it is  well known  that nonlinear methods of approximation can often perform dramatically better than
  linear methods,  there are still   questions on how to measure  the optimal performance possible for such methods.   
  This paper studies nonlinear methods of approximation that are compatible with numerical implementation in 
  that they are required to be  numerically stable.
      A measure of optimal performance, called {\em stable manifold widths},  for  approximating a model class $K$ in a 
      Banach space $X$ by stable  manifold methods is introduced. Fundamental inequalities  between these stable manifold 
      widths  and the entropy of $K$ are established.   The effects of requiring stability in the settings of deep learning and 
      compressed sensing are discussed.
      \end{abstract}


 \section{Introduction}  Nonlinear methods  are now used   in  many areas of numerical analysis, signal/image processing, and statistical learning.  
 While their improvement of error reduction when compared to linear methods is well established,  the intrinsic limitations of such  methods have not been given, 
 at least for  what numerical analysts would consider as acceptable algorithms.

   Several notions of widths have been introduced to quantify optimal performance of  nonlinear  approximation methods.    
   Historically, the first of these was the Alexandroff width described in \cite{A}.   Subsequently, alternate descriptions of  widths were given in \cite{DHM}.  
   We refer the reader to \cite{DKLT}, where  a summary of different nonlinear widths and their relations to one another is discussed.  
  
  While these 
 notions of nonlinear widths were shown to monitor certain approximation methods such as wavelet compression, they did not provide a realistic estimate for  the optimal performance of nonlinear methods in the context of numerical computation.   The key ingredient missing in these notions of widths was stability.   Stability is   essential  in numerical computation and  should be included in formulations  of the best possible performance by numerical methods.
 
   In this paper, we modify the definition of nonlinear widths to include stability.  In this way, we provide a more realistic benchmark  for  the optimal performance of numerical algorithms whose ultimate goal is to recover an underlying function.  Such algorithms are the cornerstone of numerical methods for solving operator equations, statistical methods in regression and classification, and in compressing and encoding signals and images.  It turns out that these new notions of widths  have considerable interplay with various results in functional analysis, including the bounded approximation property and the extension of Lipschitz  mappings.

 The canonical setting in approximation theory is that we are given a Banach space $X$ equipped with a norm $\|\cdot\|_X$ and we wish to approximate the elements of $X$  with error measured in this norm by simpler, less complex elements such as polynomials, splines, rational functions, neural networks,  etc.  
 The quality of this approximation is  a critical element in the design and analysis of numerical methods.   Any numerical method for computing functions is built on some form of approximation and hence the optimal performance of the numerical method   is no better than the optimal performance of the approximation method.   Note, however,
  that it may not be easy to actually design a numerical method  in a given applicative context that achieves 
 this optimal performance.   For example,   one may not be given a complete access to the target function. This is the case when we are only given limited data about the target function, as  it occurs in statistical learning and  in the theory of optimal recovery.
 
 In analyzing the performance of approximation/numerical methods, we typically examine their performance on model classes $K\subset X$,
 i.e., on  compact subsets $K$ of $X$.    The model class $K$ summarizes what we know about the target function.  For example, when numerically solving a partial differential equation (PDE), $K$ is typically provided by a regularity theorem for the PDE. In the case of signal processing, $K$ summarizes what is known or assumed about the underlying signal such as bandlimits  in the frequency domain or sparsity.

 The concept of widths was introduced to quantify the best possible performance of approximation methods on a given model class $K$.  The best known among these   widths is  the Kolmogorov width, which was introduced to quantify the best possible approximation using linear spaces.   If $X_n\subset X$ is a linear subspace of $X$ of finite dimension $n$,
 then its performance in approximating the elements of the model class $K$ is given by the {\it worst case error}
 \be
 \label{classerror}
 E(K,X_n)_X := \sup_{f\in K} \dist(f,X_n)_X.
 \ee
The value of $n$ describes the complexity of the approximation or numerical method using the space $X_n$.  If we fix the value of $n\ge 0$, the Kolmogorov $n$-width of $K$
 is defined as
 \be
 \label{Kwidth}
  d_0(K)_X=\sup_{f\in K} \|f\|_X, \quad d_n(K)_X:=\inf_{\dim(Y)=n} E(K,Y) _X, \quad  n\geq 1.
 \ee
 It  tells us the optimal performance possible on the model class $K$ using linear spaces of dimension $n$ for the approximation. Of course, it does not tell us how to select
 a (near) optimal space $Y$ of dimension $n$ for this purpose. 

 For classical model classes such as a finite ball in smoothness spaces like the Lipschitz, Sobolev, or Besov spaces,  the Kolmogorov widths are known asymptotically.  Furthermore, it is often known  that  specific linear spaces of dimension $n$ such  as polynomials,
 splines on uniform partition, etc.,   achieve this  (near) optimal performance (at least within reasonable constants).  This can then be used to show that certain numerical methods, such as spectral methods or finite element methods are also (near) optimal among all possible choices of numerical methods built on using linear spaces of dimension $n$ for the approximation.
 
 Let us note that in the definition of Kolmogorov width,  we are not requiring that   the mapping 
 which sends $f\in K$ into an approximation to $f$  is a linear map.    
  There is a concept of {\it linear width} which requires the linearity of the approximation map.  
Namely, given $n\ge 0$ and a model class $K\subset X$, its {\it linear width}   $d_n^L(K)_X$ is defined as
\be
   \label {linwidth}
     d_0^L(K)_X=\sup_{f\in K} \|f\|_X, \quad    
 d_n^L(K)_X:=\inf_{L\in \cL_n} \sup_{f\in K} \|f-L(f)\|_X, \quad  n\geq 1,
\ee
    where the infimum is taken over the class  $\cL_n$ of all  continuous linear maps from $X$ into itself with rank at most $n$.   
 The asymptotic decay of linear widths for classical smoothness classes are known.
    We refer the reader to the book of Pinkus \cite{PP} for the fundamental results for Kolmogorov and linear widths. 
 
 There is a general lower bound on the decay of   the Kolmogorov width that was given by Carl in \cite{C}.  Given $n\ge 0$, we  define the {\it entropy number} $\e_n(K) _X$ to be the infimum of all $\e>0$ for which $2^n$ balls of radius $\e$ cover $K$.  Then, Carl proved that for each $r>0$, there is a constant $C_r$ such that whenever 
 $\sup_{ m\ge 0} (m+1)^rd_m(K)_X$ is finite, then
 \be
 \label{carl}
 \e_n(K) _X\le C_r(n+1)^{-r}\sup_{ m\ge 0} (m+1)^rd_m(K)_X.
 \ee
 Thus,  for polynomial   decay rates for approximation  of the elements of $K$ by $n$ dimensional linear spaces, this decay rate  cannot be better than that of the entropy numbers of $K$.  For many  standard model classes $K$, such as finite balls in Sobolev and Besov spaces, the decay rate of $d_n(K)_X$ is
 much worse than $\e_n(K)_X$.  
 
During the decade of the 1970's, it was recognized that the performance of  approximation and numerical methods could be significantly enhanced if one uses certain nonlinear methods of approximation in place of the linear spaces $X_n$.  For example, there was the emergence of adaptive finite element methods in numerical PDEs, the   sparse  approximation from a dictionary in signal processing, and various nonlinear methods for learning.   These  new numerical methods  can be viewed as replacing 
in the construction of the numerical algorithm
the linear space $X_n$ by a nonlinear manifold $\cM_n$  depending on $n$
parameters.  For example, in place of using piecewise linear approximation on a fixed partition with $n$ cells, one would use piecewise linear approximation on a partition
of $n$ cells  which would be allowed to vary with the target function.  Adaptive finite element methods (AFEM) are a primary example of such nonlinear algorithms.
Another relevant example of nonlinear approximation,  which is  of much interest these days, are neural networks.
The parameters of the neural network are chosen depending on the target function (or the available information about the target function given through data observations) and hence is a nonlinear procedure.   The outputs of neural networks with fixed architecture form a nonlinear parametric family $\cM_n$ of functions, where $n$ is the number of parameters.  

When analyzing the performance of numerical algorithms built on some form of approximation (linear or nonlinear), an important new ingredient emerges, namely, the notion of stability.  Stability means that when the input (the information about the target function) is entered into the algorithm, the performance of the algorithm is not severely affected by small inaccuracies.  Moreover, the algorithm should not be severely effected by small inaccuracies in computation since such inaccuracies are inevitable. Having this in mind, we are interested in the following  fundamental question in numerical analysis:
\vskip .1in

\noindent
{\bf Question:}  {\it  Given a numerical task  on a model class $K$,  is there a best stable numerical algorithm for this task  and accordingly,  is there
an optimal rate-distortion performance which incorporates the notion of stability?}
\vskip .1in

\noindent
 In this context,  to formulate the notion of best,  we need a precise definition of what are admissible numerical algorithms.  We would like a notion 
  that is built on nonlinear methods of approximation   and  also respects the requirement of numerical stability.   In this paper, we take the view that nonlinear methods of approximation depending on $n$ parameters
 are  built on two mappings. 
\begin{itemize}
\item A mapping $a=a_n:X\to \R^n$,  which when given $f\in X$ chooses $n$ parameters 
 $a(f) \in \R^n$ to represent $f$. Here, when $n=0$, we take $\R^0:=\{0\}$. 

\item A mapping $M=M_n:\R^n\to X$ which maps a   vector $y\in\R^n$ back into $X$ and is used to build the approximation of $f$. 
The set 
$$
\cM_n:=\{ M_n(y):\ y\in\R^n\} \subset X
$$ 
is viewed as a parametric manifold.    
\end{itemize}

Given $f\in X$, we approximate $f$ by $A(f)=M\circ a(f):=M(a(f))$.     The error for approximating $f\in X$ is then given by
 \be 
 \nonumber
 E_{a,M}(f):=\|f-M(a(f))\|_X,
 \ee
 and the approximation error on a model class $K\subset X$ is
 \be 
 \nonumber
 E_{a,M}(K)_X:=\sup_{f\in K} E_{a,M}(f).
 \ee
 
A significant question is what conditions should be placed on  the mappings $a,M$.   If no conditions at all are placed on these mappings, we would allow  discontinuous or non-measurable mappings that have no stability and would not be useful in a numerical context.   This observation  led to requiring that  both mappings $a,M$  at least  be continuous and motivated the  definition of the  {\it manifold width  $\delta_n(K)_X$}, see \cite{DHM,DKLT},
 \be
 \label{manwidth}
 \delta_n(K)_X:=\inf_{ a,M} E_{ a,M}(K)_X,
 \ee
 where the infimum is taken over all  mappings $a:K\to \R^n$ and $M:\R^n\to X$ with  $a$ continuous on $K$ and $M$ continuous on $\R^n$.  A comparison between manifold widths and other types of nonlinear widths was given in \cite{DKLT}.   
 
Note that in numerical applications one faces the following two inaccuracies in algorithms:
 \noindent
 \begin{enumerate}
\item  In place of inputting $f$ into the algorithm, one rather inputs a noisy discretization of $f$ which can be viewed 
 as a perturbation of $f$.  So one would like to have the property that when $\|f-g\|_X$ is small
 then the algorithm outputs $M\circ a(f)$ and $M\circ a(g)$ are close to one another.  A standard quantification of this is to require that the   mapping $A:=M\circ a$ is a Lipschitz  mapping.
\item   In the numerical implementation of the algorithm the parameters $a(f)$ are not  computed exactly and so one would like 
to have the property that if $a,b\in \R^n$ are close to one another then $M(a)$ and $M(b)$ are likewise close.
 Again, the usual quantification of this in numerical implementation is that the mapping $M:\R^n\to X$ is a Lipschitz map.  This property requires the specification of a norm on $\R^n$ which is controlling the size of the perturbation of $a$.
 \end{enumerate}
  One simple way to guarantee that  these two properties  hold is to require that the  two mappings $a,M$ are themselves Lipschitz.  
  Note that this requirement implies (i) and (ii) but is indeed stronger.   We shall come back to this point 
   later in the paper.
 At present, this motivates us to introduce the following {\it stable manifold width}.  We  fix a constant $\gamma\geq 1$ and consider 
 mappings $a$ and $M$  that are {\it $\gamma$ Lipschitz continuous}   on their domains with respect to a  norm $\|\cdot\|_{Y}$ on $\R^n$, that is
 \be
 \label{Lipschitz}
 \|a(f)-a(g)\|_Y\leq\gamma \|f-g\|_X, \quad {\cnew \rm and} \quad \|M(x)-M(y)\|_X\leq \gamma \|x-y\|_{Y}, \quad x,y\in \R^n.
 \ee
 Then,  the {\it stable manifold width} $\delta_{n,\gamma}^*(K)_X$ of the compact set $K\subset X$
  is defined as
 \be
 \label{stablewidth}
 \delta_{n,\gamma}^*(K)_X:= \inf_{a,M,\|\cdot\|_{Y}} E_{a,M}(K)_X,
 \ee
 where now the infimum is taken  over all maps    $a:K\to (\R^n,\|\cdot\|_{Y})$, $M:(\R^n,\|\cdot\|_{Y})\to X$, and norms $\|\cdot\|_{Y}$ on $\R^n$, where $a,M$ 
are $\gamma$ Lipschitz.

\begin{remark} 
Note that a rescaling $\t a(f)=c a(f)$ and $\t M(x)=M(c^{-1}x)$ leaves $E_{a,M}(K)_X$ unchanged. Therefore, if 
$a$ is Lipschitz with constant $\lambda_1$ and $M$ is Lipschitz 
with constant $\lambda_2$ we can rescale them    to satisfy our definition with 
constant $\sqrt{\lambda_1 \lambda_2}$. 
We choose the above version of the definition for simplicity of notation. 
\end{remark}

Throughout the paper, we use the standard notation
\be
\ell_p^n:=(\R^n,\|\cdot\|_{\ell_p})
\ee
for the space $\R^n$ equipped with the $\ell_p$ norm,
and use $\|\cdot\|_{\ell_p^n}$ when we need to stress the dependence on $n$, or simply $\|\cdot\|_{\ell_p}$
when there is no ambiguity, for the corresponding
$\ell_p$ norm.

 The stable manifold width defined above gives a benchmark for accuracy which no  Lipschitz stable numerical algorithm can exceed when numerically recovering the model
 class $K$.  Note, however, that whether there is a numerical procedure  that can achieve this accuracy depends in part on what access is available to the 
 target functions from $K$.
 In typical numerical settings, one may not have full access to $f$ and this would restrict the possible  performance  of a numerical procedure.  For example,
 if we are only given partial information in the form of data about $f$, then performance will be limited by the quality of that data.

  The majority of this paper is   a  study   of this  stable manifold width.  We begin in the next section by discussing some of its fundamental properties.
  It turns out that some of these properties are closely connected to  classical  concepts in the theory of Banach spaces.  For example, we prove in Theorem \ref{T:tozero}
  that a separable Banach space $X$ has the property that 
   $\bar \delta_{n,\gamma}(K)_X\rightarrow 0$, $n\rightarrow\infty$, for every compact set $K\subset X$ if and only if  $X$ has the $\gamma^2$-bounded  approximation property. Here,
 $\bar \delta_{n,\gamma}(K)_X$ is a modified stable manifold width, defined the same way as  $\delta_{n,\gamma}^*(K)_X$, 
 with the only difference being  that  the infimum is taken over all $a:X\to \R^n$ defined on the whole space $X$ (rather than only on $K$) which are $\gamma$ Lipschitz.
  
  The next part of this paper seeks comparison of stable manifold widths of a compact set $K\subset X$ with its entropy numbers.   In \S \ref{sec:bounds},  we show that
   for a general Banach space $X$  stable manifold widths $\delta^*_{n,\gamma}(K)_X$ essentially cannot go to zero  faster than the entropy numbers of $K$.  Namely, we show   that for any $r>0$, we have
  \be
 \label{Carl}
  \varepsilon_n(K)_X\le C(r,\gamma)  (n+1)^{-r}\sup_{m\ge 0}(m+1)^r\delta^*_{m,\gamma}(K)_X, \quad n\ge 0.
 \ee
 Inequalities of this type are called {Carl's type inequalities} since such inequalities were first proved for Kolmogorov widths by Carl \cite{C}.  This inequality says that if 
 $\delta^*_{ n,\gamma}(K)_X$ tends to zero like $n^{-r}$ as $n$ tends to infinity,   then the entropy numbers must at least  do the same.  The significance of Carl's inequality is that in practice it is usually much easier to estimate the entropy numbers of a compact set $K$ than it is to   compute its widths.   In fact, the entropy numbers of all classical Sobolev and Besov finite balls in an $L_p$ space (or Sobolev space) are known.
Note that the assumption of stability is key here since we 
 show that less restrictive forms of nonlinear widths, for example the manifold widths, do not
 satisfy a Carl's inequality.

 While, the inequality \eref{Carl} is significant, one might speculate that in general $\e_n(K) _X$ may go to zero much faster than 
 $\delta^*_{ n,\gamma}(K){_X}$.   In \S \ref{sec:uboundsH},  we show that  when $X$ is a Hilbert space $H$, for any compact set $K \subset H$, we have
 \be
 \label{compare1}
 \delta^*_{26n,2}(K)_H\le 3 \e_n(K)_H,\quad n\ge 1.
 \ee
 We prove \eref{compare1}  by exploiting well know results from functional analysis (the Johnson-Lindenstrauss embedding lemma together with the existence of extensions of Lipschitz mappings).  When combined with the Carl's inequalities this shows that $\delta^*_{n,\gamma}(K)_H$ and $\e_n(K)_H$ behave the same when the approximation takes place in a Hilbert space $H$.
 Thus, the entropy numbers of a compact set provide a benchmark for the best possible performance of numerical recovery algorithms  in this case.
 
 A central question  (not completely answered in this paper) is what are the best  comparisons like  \eref{compare1} that  hold for a general Banach space $X$?  In section \S\ref{sec:Banach},  we prove   some first  results of the form \eref{compare1} for more general Banach spaces. Our results show some loss over \eref{compare1} when moving from a Hilbert space to a general Banach space in the sense that  the constant $3$ is now replaced by $C_0n^\alpha$,
  where $\alpha$ depends on the particular Banach space.  This topic seems to be intimately connected with the problem of extension of Lipschitz maps 
  defined on a subset $S$ of $X$ to all of $X$.
  
 From the viewpoint of approximation theory and numerical analysis, it is also of interest how classical 
  nonlinear approximation procedures comply with the stability properties proposed in this paper.  
  This is discussed in \S 6 for compressed sensing and neural network approximation. Another relevant issue is
to determine the asymptotic behavior of  $\delta^*_{ n,\gamma}(K)_X$ for classical smoothness classes $K$ used in
 numerical analysis, for example when $K$ is the unit ball of a Sobolev or Besov space. 
 For now, only in the case $X=L_2$  is there a satisfactory understanding of this behavior.



 \section{Properties of stable  manifold widths}
 
 In this section, we derive properties of the stable manifold width and discuss its relations  with  certain concepts
 in the theory of Banach spaces such as  the bounded approximation property.

 \subsection{On the definition of $\delta^*_{ n,\gamma}(K)_X$}
  Let us begin by making some comments on the definition of $\delta_{n,\gamma}^*(K){ _X}$  presented in \eref{stablewidth}.  
  In this definition, we assumed   that the mappings $a$ were Lipschitz only  on $K$.   We could have imposed the stronger condition that $a$ is defined and Lipschitz
 on all of $X$.  Since this concept is sometimes useful, we define the {\it modified stable  manifold width}
 \be
 \label{msw}
  \bar\delta_{n,\gamma}(K)_X:=\inf_{a,M,\|\cdot\|_{Y}} \sup_{f\in K} \|f-M(a(f))\|_X,
 \ee
 with the infimum now taken over all   norms $\|\cdot\|_{Y}$ on $\R^n$ and mappings $a:X\to (\R^n,\|\cdot\|_{Y})$ and $M:(\R^n,\|\cdot\|_{Y})\to X$ 
which are $\gamma$ Lipschitz. 
 Obviously, we have
 \be
 \label{compare2}
 \delta_{n,\gamma}^*(K)_X\le \bar \delta_{n,\gamma}(K)_X,\quad  n\geq 0.
 \ee
  On the other hand,  in the case of a Hilbert space $H$, the following lemma holds.
 
 \begin{lemma}
For  $K\subset H$ a  compact convex subset of the  Hilbert space $H$ we have
$$
\delta_{n,\gamma}^*(K)_H=\bar \delta_{n,\gamma}(K)_H, \quad n\geq 0.
$$ 
 \end{lemma}
 
 \noindent
{\bf Proof:}  Having in mind \eref{compare2}, we only need to show that 
$\delta_{n,\gamma}^*(K)_H\geq  \bar \delta_{n,\gamma}(K)_H$. Let us fix $n\geq 0$ and let
$$
a:K\to (\R^n,\|\cdot\|_{Y})
$$
be any $\gamma$ Lipschitz map and let us consider the 
metric projection $P_K:H\to K$ of $H$ onto $K$,
$$
P_K(f):=\displaystyle{\argmin_{g\in K}\|g-f\|_H}.
$$
Note that $P_K$ is $1$ Lipschitz map.  Therefore, $a$ can be extended to the $\gamma$ Lipschitz map
 $$
\t a:= a\circ P_K:H\to  (\R^n,\|\cdot\|_{Y})
 $$
defined on $H$, and we find that $E_{\t a,M}(K)_X=E_{a,M}(K)_X$ for any reconstruction map $M$. Thus,  $\delta_{n,\gamma}^*(K)_H\geq  \bar \delta_{n,\gamma}(K)_H$, and the proof is completed.
\hfill $\Box$

\begin{remark}
\label{r1}
The above approach relies 
on properties of metric projections, see \cite {AN}, and can be used to show  intrinsic relations between 
$\delta_{n,\gamma}^*(K)_X$ and $\bar\delta_{n,\gamma}(K)_X$ for certain  compact subsets $K\subset X$ of a general Banach space $X$. 
\end{remark}

   \begin{remark}
   \label{L:newdef}
  In the definition of $\delta^*_{ n,\gamma}(K)_X$ and $\bar \delta_{ n,\gamma}(K)_X$,  the space $(\R^n, \|\cdot\|_{ Y})$ 
  can be replaced by any normed space $(X_n,\|\cdot\|_{X_n})$ of dimension $n$.  That is,  for example, in the case of $\delta^*_{n,\gamma}(K)_X$,
  \be
  \label{modifieddelta}
  \delta^*_{n,\gamma}(K)_X= \inf_{a,M,X_n}\sup_{f\in K} \|f-M(a(f))\|_X,
  \ee
  where now the infimum is taken over all  normed spaces $X_n$ of dimension $n$ with norm $\|\cdot\|_{X_n}$ and 
  all $\gamma$ Lipschitz  maps $a:X\to (X_n,\|\cdot\|_{X_n})$ and 
  $M:(X_n,\|\cdot\|_{X_n})\to X$. Indeed, consider any basis $(\phi_1,\dots,\phi_n)$ of $X_n$. The associated coordinate map 
$\kappa: X_n\to \R^n$ defined by {\cnew $\kappa(g)=(x_1,\dots,x_n)=x$} for $g=\sum_{i=1}^n x_i\phi_i$ is
an isometry when $\R^n$ is equiped with the norm $\|x\|_Y:=\|g\|_{X_n}$. For this norm, the maps
$\t a=\kappa \circ a:X\to \R^n$ and $\t M=M \circ \kappa^{-1}:\R^n\to X$ have the same Lipschitz constants as $a:X\to X_n$ and 
$M: X_n\to X$,
which shows the equivalence between the two definitions.
\end{remark}

  \subsection {When does $\bar \delta_{n,\gamma} (K)_X$ tend to zero as $n\to \infty$?}
 \label{SS:tozero}
 
 
 We turn next to the question of whether $\bar \delta_{n,\gamma}(K)_X$ tends to zero for all compact sets $K\subset X$.  In order to orient this discussion, we first recall results of this type
 for other widths and for other closely related concepts in the theory of Banach spaces.

  Let $X$ be a separable Banach space.  While the Kolmogorov widths $d_n(K)_X$ tend to zero as $n\to \infty$ for each compact set $K\subset X$, notice that this definition of widths says nothing about how the approximants to a given $f\in K$ are constructed.
 In the definition of  the  
  linear widths  $d^L_n(K)_X$, see \eref{linwidth},  
 it is required that the approximants to $f$ are constructed by finite rank   continuous linear mappings. 
  In this case, it is known that a necessary and sufficient condition that these widths tend to zero is that $X$ has the {\it  approximation property}, i.e. 
  for each compact subset $K\subset X$,
  there is a sequence of 
 bounded linear  operators $T_n$ of finite
 rank  at most $n$ such that 
 \be
 \label{ap}
 \sup_{f\in K} \|f-T_n(f)\|_X\to 0,\quad n\to \infty.
 \ee
 In the definition of  approximation property,  the norms of the  operators $T_n$ are allowed to grow with $n$.  
 A second concept of  $\gamma$-{\it bounded  approximation property } requires  in addition that there is a $\gamma\geq 1$
 such that  the operator norm bound $\|T_n\|\le \gamma$ holds for the operators in \eref{ap}.

 The main result of this section is the  following theorem which  characterizes the Banach spaces $X$ for which every 
 compact  subset $K\subset X$ has the property  $\bar \delta_{n,\gamma} (K)_X\to 0$ as $n\to\infty$.

\begin{theorem}
\label{T:tozero}
Let $X$ be a separable Banach space and  $\gamma \geq 1$. The following two statements are equivalent:
\vskip .05in
\noindent
\hskip 0.2in{\rm (i)}  $\bar\delta_{n,\gamma}(K)_X\rightarrow 0$ as $n\to \infty$ for every compact set  $K\subset X$.
\vskip .05in
\noindent
\hskip 0.2in{\rm (ii)}    $X$ has the $\gamma^2$-bounded  approximation property.
\end{theorem} 
Before going further, we state a lemma that we use in the proof of the above theorem. The proof of the lemma is given  after  the proof of the theorem.
\begin{lemma}
\label{L:modified}
Let  $\|\cdot\|_Y$ be a norm on $\R^n$, $n\geq 1$, and  $X$ be any separable Banach space.    If
 $M:\R^n\to X$ is  a $\gamma$ Lipschitz  mapping,  then for any bounded set $S\subset  \R^n$, and any $\e>0$, 
  there exists a map $\o M:  \R^n\to X$,  $\o M=\o M(S,\e)$, 
with the following properties:
\vskip .05in
\noindent
\hskip 0.2in{\rm (i)} $\o M$ is  Lipschitz with constant $\gamma$.
\vskip .05in
\noindent
\hskip 0.2in{\rm (ii)}   $\o M$ has finite rank, that is    $\o M(\R^n)$ is a subset of a finite dimensional subspace of $X$.

\noindent
\hskip 0.2in{\rm (iii)}  $\o M$ approximates $M$ to accuracy $\e$ on $S$, namely
\be
\label{perturbo}
\|M-\o M\|_{L_\infty(S,X)}:=\max_{x\in S} \|M(x)-\o M(x)\|_X \leq \e.
\ee
 \end{lemma}
 \vskip .1in
 \noindent
{\bf Proof of Theorem \ref{T:tozero}:}   First, we show that (ii) implies (i).   If $X$ has the $\gamma^2$-bounded approximation property, then given any compact set $K\subset X$, there is a sequence of operators $\{T_n\}$, $n\ge 1$, $T_n:X\to X_n$ with $X_n$ of dimension at most $n$,
 with operator norms $\|T_n\|\le\gamma^2$, and 
\be
\label{T01}
\sup_{f\in K} \|f-T_n(f) \|_X\to  0,\quad n\to \infty.
\ee
Consider the mappings 
$$
a:= \gamma^{-1} T_n:X\to X_n, \quad M:= \gamma  Id:X_n\to X_n\subset X.
$$
Each of these mappings is Lipschitz with Lipschitz constant at most $\gamma$ and $M\circ a=T_n$. By virtue of \eref{T01} and 
Remark \ref{L:newdef}, we have that $\bar \delta_{n,\gamma}(K)_X\rightarrow 0$ as $n\rightarrow 0$.

 Next, we show that (i) implies (ii).  Suppose that (i) of Theorem \ref{T:tozero} holds and $K$ is any compact set in $X$.   
From the definition of $\bar \delta_{n,\gamma}(K)_X$, there exist  $\gamma$ Lipschitz mappings
$$
a_n:X\to  \R^n, \quad M_n: \R^n\to X,
$$ 
with some norm $\|\cdot\|_{Y_n}$ on $\R^n$ and 
%
\be
\nonumber
\sup_{f\in K}\|f-M_n\circ a_n(f)\|_X\to 0,\quad n\to\infty.
\ee
We take  $\varepsilon =1/n$ in Lemma \ref{L:modified}  and let  $\o M_n$ be the modified mapping  for $M_n$  guaranteed by   the lemma with the set $S$ being $a_n(K)$.  
Then the mapping $T_n:X\to X$ defined by 
$$
T_n:=  \o M_n\circ a_n
$$
is  $\gamma^2$ Lipschitz and  has a finite rank.  Moreover, 
since for every $f\in K$,
$$
\|f-T_n(f)\|_X\leq \|f-M_n\circ a_n(f)\|_X+\|M_n\circ a_n(f)-\o M_n\circ a_n(f)\|_X\leq \|f-M_n\circ a_n(f)\|_X+1/n,
$$
one has
\be
\nonumber
\sup_{f\in K}\|f-T_n(f)\|_X\to 0,\quad n\to\infty.
\ee
To complete the proof, we  use  Theorem 5.3 from \cite{GK}, see also the discussions in \cite{G, G1},  to conclude that  $X$  has the $\gamma^2$-bounded approximation property.
\hfill $\Box$
\nl

 We now proceed with the proof of the lemma.
 \vskip .1in

\noindent
{\bf Proof of Lemma \ref{L:modified}:}
We fix the value of  $n\geq 1$  and a norm $\|\cdot \|_Y$ on $\R^n$. 
We will  prove the apparently weaker statement that for any
$\e,\delta>0$, there exists a $(\gamma+\delta)$ Lipschitz map $\wt M: \R^n\to X$ 
with finite rank such that
\be
\|M-\wt M\|_{L_\infty(S,X)}:=\max_{x\in S} \|M(x)-\wt M(x)\|_X \leq \e.
\label{perturb}
\ee
Once we construct $\wt M$,  we   obtain
the claimed statement by taking
\be
\nonumber
\o M=\frac {\gamma}{\gamma+\delta} \wt M.
\ee
 Clearly,  $\o M$ will satisfy (i), (ii), and (iii), since 
\begin{eqnarray}
\nonumber
\|M-\o M\|_{L_\infty(S,X)}&\leq &\|M-\wt M\|_{L_\infty(S,X)} +\|\wt M-\o M\|_{L_\infty(S,X)}\\ \nonumber
&\leq& \e+\frac{\delta}{\gamma+\delta}\max_{x\in S} \|\wt M(x)\|_X< \e+\frac{\delta}{\gamma}\max_{x\in S} \|\wt M(x)\|_X,
\nonumber
\end{eqnarray}
where $\delta$ and $\e$ are arbitrarily small and $S\subset \R^n$ is bounded.

The construction of $\wt M$ from $M$ proceeds in $3$ steps, where  one of the main issues is to keep  control of the Lipschitz constants.
\vskip .1in
\noindent
{\bf Step 1:}  Let us fix $\delta>0$.
In this step, we  construct  a map $M_1$
that agrees with $M$ on $S$, takes
the constant value $M(0)$ outside of a larger set  that contains $S$, and is $(\gamma+\delta/2)$ Lipschitz.  
We take $R_1>0$ sufficiently large such that 
$S$ is contained in the ball of radius $R_1$ with respect to the $\|\cdot\|_Y$ norm, that is,
\be
\nonumber
x\in S\implies \|x\|_Y  < R_1.
\ee
For  $\lambda>0$,  we then define  the continuous piecewise linear  function $ \phi_\lambda:\R^+\to \R$ by
\begin{equation}
\nonumber
{ \phi_\lambda}(t)=\begin{cases}
1,\quad \quad \quad \quad \quad \quad \,\,\,0\leq t\leq R_1,\\
1-\lambda(t-R_1), \quad R_1\le t\le R_1+1/\lambda,\\
0, \quad \quad \quad \quad \quad \quad \,\, \,t\geq R_1+1/\lambda.
\end{cases}
\end{equation}
 Clearly, $\phi_\lambda$ is $\lambda$ Lipschitz function and $0\leq \phi_\lambda(t)\leq 1$ for all $t\geq 0$. Next, we define the function 
 $\Phi_\lambda:\R^n\to \R^n$ by
 \be
 \nonumber
  \Phi_\lambda(x):= \phi_\lambda(\|x\|_Y)x =
\begin{cases}
x, \quad \quad \quad \quad \quad \quad \quad \quad \,\,\,\quad \|x\|_Y\leq R_1,\\
(1-\lambda(\|x\|_Y-R_1))x, \quad R_1\le \|x\|_Y\le  R_1+1/\lambda,\\
0,\quad \quad \quad \quad \quad \quad \quad \quad \,\,\, \quad \|x\|_Y\geq R_1+1/\lambda,
\end{cases}
 \ee
 and thus $\Phi_\lambda(x)=x$ for $x\in S$.
Let us check the Lipschitz property of $\Phi_\lambda$.  

First,  for $x,y$ contained in the ball  $B$ of radius $R_1+1/\lambda$ with respect to the $\|\cdot\|_Y$ norm, we have
 \be
\nonumber
  \Phi_\lambda(x)-\Phi_\lambda(y)=    \left (\phi_\lambda(\|x\|_Y)- \phi_\lambda(\|y\|_Y)\right )x+  \phi_\lambda(\|y\|_Y)(x - y),
\ee
and  thus
  \begin{eqnarray}
 \label{write1}
\| \Phi_\lambda(x)-\Phi_\lambda(y)\|_Y & \le &\| x\|_Y\big|\phi_\lambda(\|x\|_Y) -\phi_\lambda(\|y\|_Y)\big |+  \phi_\lambda(\|y\|_Y)\| x-y\|_Y
\nonumber\\
&\le & \lambda \|x\|_Y |\|x\|_Y-\|y\|_Y|+ \phi_\lambda(\|y\|_Y) \| x-y\|_Y
\nonumber \\
&{\cnew \leq}& (\lambda\|x\|_Y+\phi_\lambda(\|x\|_Y))\|x-y\|_Y \nonumber\\
&\leq& (1+\lambda R_1) \|x-y\|_Y,\quad  x,y\in B.
\end{eqnarray}

Next, for $x,y\in \R^n$ such that $\|x\|_Y\geq R_1+1/\lambda$ and  $\|y\|_Y\geq R_1+1/\lambda$
\be
\label{pop}
 \Phi_\lambda(x)-\Phi_\lambda(y)=0.
\ee
Lastly,  if  $\|x\|_Y\leq R_1+1/\lambda$ and  $\|y\|_Y>R_1+1/\lambda$, we consider the  point $x^*:=x+s^*(y-x)$, $s^*\in[0,1]$ of the intersection of the line segment connecting 
$x$ and $y$ and the sphere with radius $R_1+1/\lambda$. We have $\Phi_\lambda(y)=\Phi_\lambda(x^*)=0$, and thus it follows from \eref{write1} that
\begin{eqnarray}
\label{pop1}
\| \Phi_\lambda(x)-\Phi_\lambda(y)\|_Y&=& \|\Phi_\lambda(x)-\Phi_\lambda(x^*)\|_Y\leq 
(1+\lambda R_1)\|x-x^*\|_Y\\ \nonumber
&=& (1+\lambda R_1)s^*\|x-y\|_Y\leq (1+\lambda R_1)\|x-y\|_Y.
\end{eqnarray}

From   \eref{write1}, \eref{pop}, and \eref{pop1}, we conclude  that $\Phi_\lambda$ is a $(1+\lambda R_1)$ Lipschitz function.
 We can make the Lipschitz constant  $(1+\lambda R_1)$   as close to one as we wish by taking $\lambda$ small.
 Therefore,  choosing $\lambda$ sufficiently small, we have that  the function
\be
\nonumber
M_1:=M\circ \Phi_\lambda,
\ee
is $(\gamma+\delta/2)$ Lipschitz, agrees with $M$ over $S$
and has constant value $M(0)$  on  the set 
$$
\{x\in \R^n:\,\|x\|_Y\geq R_1+1/\lambda\}.
$$ 
 By equivalence of norms {\cnew on $\R^n$}, we conclude that $M_1$
has value $M(0)$ outside an $\ell_\infty$ cube $ [-R_2,R_2]^n$, 
{\cnew with $R_2=R_2(\lambda,n)$}

\vskip .1in
\noindent
{\bf Step 2:} In the second step,  we approximate $M_1$ by a function 
$M_2$ obtained by regularization,  see \cite{HNVW}. We consider a standard mollifier
 \be
 \nonumber
\vp_m(x)=m^{n} \vp(m x), \quad x\in\R^n,
\ee 
where $\vp$ is a smooth positive function supported on the unit euclidean ball of $\R^n$
and such that $\int_{\R^n}\vp=1$. We then define $M_2:=\vp*M_1$, that is,
\be
\nonumber
M_2(x)=M_2(m,x):=\int_{\R^n} \vp_m(y)M_1(x-y)dy.
\ee
The function $M_2$ is smooth and equal to $M(0)$ outside of the cube 
\be
Q:=[-D,D]^n, \quad\quad {\cnew D:=R_2+\frac1m}.
\label{cube}
\ee

By taking $m$ sufficiently large, we are ensured that 
\be
\nonumber
\max_{x\in \R^n} \|M_1(x)-M_2(x)\|_X \leq \e/2,
\ee
and in particular (since $M_1$ agrees with $M$ on $S$)
\be
\max_{x\in S} \| M(x)-M_2(x)\|_X \leq \e/2,
\label{perturb1}
\ee
because
\begin{eqnarray}
\nonumber
\| M_1(x)-M_2(x)\|_X&=&\Big\|\int_{\R^n} \vp_m(y)M_1(x)dy-\int_{\R^n} \vp_m(y)M_1(x-y)\,dy\Big\|_X\\
\nonumber
&\leq &\int_{\R^n} \vp_m(y)\|M_1(x)-M_1(x-y)\|_X\,dy\\
\nonumber
&\leq & (\gamma+\delta/2)\int_{\R^n} \vp_m(y)\|y\|_Y\,dy=\frac{\gamma+\delta/2}{m}\int_{\R^n} \vp(y)\|y\|_Y\,dy.
\nonumber
\end{eqnarray}
In addition, by convexity  we find that $M_2$ is $(\gamma+\delta/2)$ Lipschitz since
\begin{eqnarray}
\label{Lip1}
\nonumber
\|M_2(x)-M_2(y)\|_X  
&=&\Big\| \int_{\R^n} \vp_m(z)(M_1(x-z)-M_1(y-z))\, dz\Big\|_X  \nonumber \\
& \leq &
\int_{\R^n} \vp_m(z) \|M_1(x-z)-M_1(y-z)\|_X \, dz   \nonumber \\
&\leq & \int_{\R^n} \vp_m(z) (\gamma+\delta/2) \|x-y\|_Y  \,dz=(\gamma +\delta/2)\|x-y\|_Y.
\nonumber
\end{eqnarray}
If we take $m$ sufficiently large then \eref{perturb1} holds and the construction of $M_2$ in this step is complete. We fix $m$  for the remainder of the proof.
Any constants $C$ given below depend only on $m$, $n$, $\delta$, and the initial function $M$.  The value of $C$ may change from line to line.
\vskip .1in

\noindent
{\bf Step 3:} In this step, we derive $\wt M$ from $M_2$ by piecewise  linear interpolation.
 We work on
the support cube $Q=[-D,D]^n$.
We recall that  $M_2$ is constant and equal to $M(0)$ outside of $Q$. We create
a simplicial mesh of  {\cnew $2Q$} by subdividing it into subcubes $Q_k$
of equal side length  $h=2D/N$, and  then using the Kuhn simplicial decomposition of each of these subcubes
into  $n!$ simplices, see \cite{K}.  The set of vertices of the cubes $Q_k$ form a mesh of discrete points in  {\cnew $2Q$}. We denote  by  
$\Lambda_h=\{x_\nu\}\subset Q\subset \R^n$
the set of these vertices that belong to $Q$ . 

We denote by $I_h$ the operator of  piecewise linear interpolation at the vertices of  $\Lambda_h$. It is usually applied to scalar valued functions but its extension to 
Banach space valued functions is immediate.  Since $M_2$ has value $M(0)$ on $\partial Q$, the same holds for $I_hM_2$ which may be written as
\be
\nonumber
\wt M (x):=I_h M_2 (x)=\sum_{x_\nu\in{\cnew \Lambda_h}}  M_2(x_\nu) N_\nu  (x),\quad x\in Q.
\ee
 Here, 
the functions $N_\nu$ are the nodal basis for piecewise linear interpolation,  that is
$N_\nu$ is a continuous piecewise linear function with  $N_\nu(x_\mu)=\delta_{\mu,\nu}$, with $\delta_{\mu,\nu}$ the Kronecker symbol
for $\mu\in \Lambda_h$. We then can extend $\wt M$  by the value $M(0)$ outside of $Q$. It follows that  
$\wt M(\R^n)$ is contained in a linear subspace of dimension  {\cnew $\#( \Lambda_h)+1$},
that is  $\wt M$ has finite rank. 
 We are now left  to show that 
 $\wt M$ is $(\gamma+\delta)$ Lipschitz and that  \eref{perturb} holds. Thus it is enough to show that:
\begin{enumerate}
\item  $(\wt M-M_2)$ is $\delta/2$ Lipschitz, since $M_2$ is $(\gamma+\delta/2)$ Lipschitz;
\item  $ \max_{x\in S}\|\wt M(x)-M_2(x)\|_X\leq \e/2$, because of \eref{perturb1}.
\end{enumerate}

In order to prove (i) and (ii), we first note  that if $U$ is the unit ball in $X^*$, we have that 
$$
\|(\wt M(x)-M_2(x))-(\wt M(y)-M_2(y))\|_X=\sup_{\ell\in U} |\ell(\wt M(x))-\ell(M_2(x))-(\ell(\wt M(y))-\ell(M_2(y)))|,
$$
and
$$
 \|\wt M(x)-M_2(x)\|_X=\sup_{\ell\in U}|\ell(\wt M(x))-\ell(M_2(x))|.
$$
 For any  $\ell\in U$, we denote by $v_\ell:\R^n\to \R$ the piecewise linear scalar valued function
\begin{eqnarray}
\nonumber
v_\ell(x):=\ell(M_2(x))=\begin{cases}
\ell(M_2(x)), \quad x\in Q,\\
\ell(M(0)), \quad x\in\R^n\setminus Q.
\end{cases}
\end{eqnarray}
Then we have 
\begin{eqnarray}
\nonumber
\ell(\wt M(x))=\begin{cases}
\ell(I_hM_2(x))=I_hv_\ell(x),\quad x\in Q,\\
\ell(M(0)), \quad \quad \quad\quad \quad \quad \quad \,\, x\in \R^n\setminus Q,
\end{cases}
\end{eqnarray}
and
\begin{eqnarray}
\nonumber
\ell(M_2(x))-\ell(\wt M(x))=\begin{cases}
v_\ell(x)-I_hv_\ell(x),\quad x\in Q,\\
0, \quad \quad \quad\quad \quad \quad \quad \,\, x\in \R^n\setminus Q.
\end{cases}
\end{eqnarray}
Note here that we have used the slight abuse of notation since the same notation $I_h$ is used for the interpolation operator applied to scalar valued functions  as well as for  Banach space valued functions.  In particular, we may extend $I_hv_\ell$ by $\ell(M(0))$ outside of $Q$.

 Therefore, to show (i) and (ii), it is enough to show that uniformly in $\ell\in U$, for $h$ sufficiently small,
 $(v_\ell-I_hv_\ell)$ is $\delta/2$ Lipschitz function on $\R^n$ and

\be
\label{ki}
\max_{x\in Q}|v_\ell(x)-I_hv_\ell(x)|\leq \e/2.
\ee
 Note that the functions $v_\ell$, $\ell \in U$,  are smooth with uniformly bounded (in $\ell$) second derivatives 
$$
|v_\ell |_{W^{2,\infty}(\R^n)}:=\max_{x\in \R^n}\max_{|\alpha|=2}|\partial_\alpha v_\ell(x)|\le C_0,.
$$
with $C_0$ a fixed constant independent of $\ell$.

Let $\cK$ be any  one of the simplices in the Kuhn simplicial decomposition of any of the $Q_k$'s. 
Then the diameter of $\cK$ is $\sqrt{n}h$ and the radius of the inscribed sphere is $\frac{h}{\sqrt{2}(n-1+\sqrt{2})}$, see  subsection 3.1.4 in \cite{K}. 
It follows from  Corollary 2 in \cite{CW} that 
$$
\max_{x\in \cK}|v_\ell(x)-I_hv_\ell(x)|\leq Ch^2|v_\ell |_{W^{2,\infty}(\cK)}, \quad 
\max_{i=1,\ldots,n}\max_{x\in \cK}|\frac{\partial}{\partial x_i}v_\ell(x)- \frac{\partial}{\partial x_i}I_hv_\ell(x)|\leq Ch|v_\ell |_{W^{2,\infty}(\cK)},
$$
 and in turn
$$
\max_{x\in \R^n}|v_\ell(x)-I_hv_\ell(x)|\leq Ch^2|v_\ell |_{W^{2,\infty}(\R^n)}, \quad 
\max_{i=1,\ldots,n}{\cnew \|\frac{\partial}{\partial x_i}v_\ell- \frac{\partial}{\partial x_i}I_hv_\ell\|_{L_\infty(\R^n)}}\leq Ch|v_\ell |_{W^{2,\infty}(\R^n)}.
$$
Thus \eref{ki} follows from the first inequality if we select $h$ small enough.
From the second of these inequalities, we find that for $x,y\in \R^n$,
\begin{eqnarray}
\nonumber
|(v_\ell(x)-I_hv_\ell(x))-(v_\ell(y)-I_hv_\ell(y))|&\le&
 C h |v_\ell|_{W^{2,\infty}(\R^n)} \|x-y\|_{\ell_1(\R^n)}
 \\ \nonumber
 &\leq& C h |v_\ell|_{W^{2,\infty}(\R^n)}\|x-y\|_Y<
\delta/2\|x-y\|_Y,
\end{eqnarray}
where we have used the fact that any two norms on $\R^n$ are equivalent and  $h$ can be made small enough. This completes the proof of the lemma.
\hfill $\Box$

\subsection{When is $\delta^*_{ n,\gamma}(K)_X=0$?}
\label{SS:deltazero}

In this section, we characterize the sets $K$ for which $\delta_{n,\gamma}^*(K)_X=0$.   We also consider a closely related question of whether 
$\delta_{ n,\gamma}^*(K)_X$ is assumed.
We will use the following lemma which is a form of Ascoli's theorem.

  \begin{lemma}
  \label{Ascoli}  Let 
  $(X,d)$ be a separable metric space and $(Y,\rho)$ be  a metric space
for which  every closed ball is compact. Let $F_n:X\to Y$ be a
sequence of $\gamma$ Lipschitz maps  for which  there exists $a\in X$
and $b\in Y$ such that $F_n(a)=b$ for $n=1,2,\dots$. Then, there exists
a subsequence $F_{n_j}$, $j\ge 1$,  which is  point-wise convergent to a function $F:X\to Y$ and
$F$ is $\gamma$ Lipschitz.  If $(X,d)$  is also compact, then the convergence is uniform.
  \end{lemma}
  
  \noindent
  {\bf Proof:} For any $f\in X$ we have
$$
\rho(F_n(f),b)=\rho(F_n(f),F_n(a))\leq \gamma d(f,a).
$$
 Let us fix a
countable dense subset $A=\{f_j\}_{j=1}^\infty\subset X$  and define
$$
B_j:=B(b,\gamma d(f_j,a))
$$ 
as the  closed  ball in $Y$ with radius $\gamma d(f_j,a)$ centered at $b$. Then the  cartesian product
$$
\cB:=B_1\times B_2\times \cdots
$$
 is a compact metric space  under  the natural
product topology. We naturally
identify each $F_n$ with an element $\hat F_n\in \cB$  whose $j$-th coordinate is
$F_n(f_j)$, that is
$$
\hat F_n:=(F_n(f_1), F_n(f_2), \ldots,F_n(f_j),\ldots )\in \cB.
$$
So, there exists a subsequence $ \hat F_{n_s}$
convergent to
 an element $ \hat F\in \cB $,  that is
 $$
 \hat F(j)=\lim_{s\to \infty} \hat F_{n_s}(j)=\lim_{s\to \infty} F_{n_s}(f_j), \quad j\geq 1.
$$

In other words,   we get a function $ F :A\to Y$, defined as
  $$
  { F}(f_j)=\lim_{s\to\infty} F_{n_s}(f_j).
  $$
   We check that
$\rho({F}(f_j),{ F} (f_i))\leq \gamma d(f_j,f_i)$ for each  $i,j=1,2,\dots$. Since $A$ is dense,   $F$ extends to a $\gamma$ Lipschitz function on $X$, $F:X\rightarrow  Y$.  Moreover, $F(f)=\lim_{s\to\infty} F_{n_s}(f)$ for every $f\in X$.  If $(X,d)$ is compact, uniform convergence is proved
by a standard argument, remarking that for any $\e>0$ we can cover $X$ by a finite number of $\e$-balls with centers $g_1,\dots,g_k\in X$,
and so
$$
\sup_{f\in X}\rho(F_{n_s}(f),F(f))\leq 2\gamma \e +\max_{i=1,\dots,n} \rho(F_{n_s}(g_i),F(g_i))\leq (2\gamma+1)\e,
$$
for $s$ large enough.
\hfill $\Box$

\begin{theorem} 
\label{T:deltazero}
Let $K\subset X$ be  a compact set in  a separable Banach space $X$.  If $\delta_{n,\gamma}^*(K)_X=0$, 
then the set  $K$ is  $\gamma$ Lipschitz
equivalent to a subset of $\R^n$. That is,   there is a norm $\|\cdot\|_{Y}$ on $\R^n$ and a    
function $F:K\to (\R^n, \|\cdot\|_{Y})$  such that $F$ is invertible and both $F$ and $F^{-1}$ are $\gamma$ Lipschitz.
    \end{theorem}
\noindent
  {\bf  Proof:} Notice that if we knew that $\delta^*_{n,\gamma}(K)_X=0$ was  assumed by maps $a,M$,  then we could simply take $F:=a$ and $F^{-1}= M|_{a(K)}$.  
So the proof consists of  a limiting argument.   
For each $k\ge 1$,  there exist a norm $\|.\|_{Y_k}$ on $\R^n$ and $\gamma$ Lipschitz maps
$a_k:K\to (\R^n,\|.\|_{Y_k})$ and $M_k:(\R^n,\|.\|_{Y_k})\to X$ such that
  \beqn 
 \label{PW.27.07}
  \lim_{k\to \infty} \sup_{f\in K}\|f-M_k(a_k(f))\|_X =\delta^*_{n,\gamma}(K)_X=0.
  \eeqn
    Let us  fix $f_0\in K$ and define 
 $$
 a_k'(f):=a_k(f)-a_k(f_0) \quad {\rm and} \quad M'_k(x):=M_k(x+a_k(f_0)).
 $$
Then, $a'_k:K\to  (\R^n,\|.\|_{Y_k})$ and $M'_k:(\R^n,\|.\|_{Y_k})\to X$ are $\gamma$ Lipschitz
 maps. Moreover,  $a'_k(f_0)=0$ and $M'_k\circ a'_k=M_k\circ a_k$ for $k=1,2,\dots$.

We  denote by $U$ the unit  ball in $\R^n$ with respect to the Euclidean norm  $\|.\|_{\ell_2^n}$   and by $U_k$ the unit ball of $\R^n$ with respect to the norm $\|\cdot\|_{Y_k}$.   From the Fritz John  theorem  (see  e.g \cite[Chapt. 3]{P}) we infer that there exist invertible linear  operators
$\Lambda_k$ on $\R^n$ such that
  $$
  U \subset\Lambda_k(U_k)\subset  \sqrt n U,
  $$
and therefore the modified norm {\cnew $\|.\|_{Z_k}$} defined as
$$
\|x\|_{Z_k}:=\|\Lambda_k^{-1}(x)\|_{Y_k}, \quad x\in \R^n,
$$
satisfies the inequality
  \beqn \label{John}
  \frac{1}{\sqrt n}\|x\|_{\ell_2^n}\leq \|x\|_{Z_k}\leq \|x\|_{\ell_2^n}.
  \eeqn
  Next, we replace $ a'_k$  and $M'_k$ by
 $$
 \tilde a_k:=\Lambda_k\circ a'_k:K \to (\R^n,\|.\|_{Z_k}), \quad \tilde M_k:=M'_k\circ \Lambda^{-1}_k:(\R^n,\|.\|_{Z_k})\to X.
$$
Note that  $M_k\circ a_k=M'_k\circ a'_k=\tilde M_k\circ \tilde a_k$, and    $\tilde a_k(f_0)=0$.  We note that   $\tilde a_k$ and $\tilde M_k$ are $\gamma$ Lipschitz with respect to the new norm $\|\cdot\|_k$. 
Indeed,
\be
\label{Lipatilde}
\|\tilde a_k(f)-\tilde a_k(g)\|_{Z_k}=\|\Lambda_k\circ a'_k(f)-\Lambda_k\circ a'_k(g)\|_{Z_k}=\|a'_k(f)-a'_k(g)\|_{Y_k}\leq \gamma \|f-g\|_X,
\ee
and
$$
\|\tilde M_k(x)-\tilde M_k(y)\|_X=\|M'_k\circ \Lambda^{-1}_k(x)-M'_k\circ \Lambda^{-1}_k(y)\|_X\leq \gamma \|\Lambda^{-1}_k(x)-\Lambda^{-1}_k(y)\|_{Y_k}=
\gamma \|x-y\|_{Z_k}.
$$

 We then extract subsequence of these mappings that converge point-wise by using  Lemma \ref{Ascoli}.  For this, we first note that from \eref{John}, we have 
 $$
 \|\tilde a_k(f)-\tilde a_k(g)\|_{\ell_2^n}\leq\sqrt{n}\|\tilde a_k(f)-\tilde a_k(g
 )\|_{Z_k}\leq \gamma\sqrt{n}\|f-g\|_X.
 $$
Hence, the sequence $\tilde a_k: K\to{\cnew\ell_2^n}$ is  a sequence of $\gamma\sqrt{n}$ Lipschitz  mappings
for which $\tilde a_k(f_0)=0$. We apply  Lemma \ref{Ascoli} to infer that,  up to a subsequence extraction, $\t a_k$ converges  point-wise  on  $K$ to a  mapping $F$. Note that $F: K\to\ell_2^n$ is also
a $\gamma\sqrt{n}$ Lipschitz map.

The remainder of the proof is to show that the function $F$ is the  mapping claimed by the theorem.
To prove it, we want first  to extract a single norm to use in place of the $\|\cdot\|_{Z_k}$.
   For this, we apply  Lemma \ref{Ascoli} to  the  subsequence  of  norms 
   $$
 \|\cdot \|_{Z_k}:{\cnew \ell_2^n}\to \R,\quad j=1,2,\ldots,
$$
viewed as   
$1$ Lipschitz
 functions,  to derive that,  up to another subsequence extraction,
 $\|\cdot\|_{Z_k}$  converges pointwise to a   
$1$ Lipschitz function from ${\cnew \ell_2^n}$ to $\R$ which we denote by $\|.\|_Y$. 
It  is easy to check that $\|\cdot\|_Y$ is  a norm on $\R^n$ and it   satisfies   
\be
\label{compareY}
\frac{1}{\sqrt n}\|x\|_{\ell_2^n}\leq \|x\|_Y\leq \|x\|_{\ell_2^n}, \quad x\in \R^n.
\ee 
We now verify the required Lipschitz properties of $F$ with respect to $\|\cdot\|_Y$.  First, we claim that
 \be
 \label{FLip}
\|F(f)-F(g)\|_Y\leq \gamma \|f-g\|_X,
\ee
namely, $F:K\to (\R^n,\|.\|_Y)$ is a $\gamma$ Lipschitz mapping.   Since $\lim_{k\to\infty} \|F(f)-\tilde a_k(f)\|_Y=0$ for all $f\in K$
because of \eref{compareY}, we prove \eref{FLip}
by showing that for any $\e>0$ and $f,g\in K$, we have
\be
\label{enoughto}
\|\tilde a_k(f)- \tilde a_k(g)\|_Y\le  \gamma \|f-g\|_X+\e,
\ee
for any sufficiently large $k$.   Now the set $S:=\{z \in \R^n:  z=\tilde a_k(f)- \tilde a_k(g),  f,g\in K, \ k\ge 1\}$ is bounded and therefore 
\be
\label{therefore1}
\sup_{z\in {  S}} |\|z\|_{Z_k}- \|z\|_Y|:= \e_k\rightarrow 0, \quad k\to \infty.
\ee
 This gives
\be
\|\tilde a_k(f)- \tilde a_k(g)\|_Y\le   \|\tilde a_k(f)- \tilde a_k(g)\|_{Z_k}+\e_k\le \gamma \|f-g\|_X+\e_k,\quad k\ge 1,
\ee
where we have used \eref{Lipatilde}.  Choosing $k$ sufficiently large we have \eref{enoughto} and in turn have proved \eref{FLip}.

  Finally, we need to check  that $F$ has an inverse on $F(K)$ which is $\gamma$ Lipschitz.   Let  $f,g\in  K$. For   every  $k$
we have  
  \begin{eqnarray}
  \nonumber
  \gamma \|F(f)-F(g)\|_{Z_k} &\geq
&\|\tilde M_k(F(f))-\tilde M_k(F(g))\|_X \geq  \|\tilde M_k(\tilde a_k(f))-\tilde M_k(\tilde a_k(g))\|_X\nonumber\\
  &-&\|\tilde M_k (F(f))-\tilde M_k(\tilde a_k (f))\|_X-\|\tilde M_k (F(g))-\tilde M_k(\tilde a_k (g))\|_X \nonumber \\
  &\geq& \|\tilde M_k(\tilde a_k(f))-\tilde M_k(\tilde a_k(g))\|_X-\gamma\|F(f)- \tilde a_k (f)\|_{Z_k}-\gamma\|F(g)- \tilde a_k (g)\|_{Z_k}.
  \nonumber
  \end{eqnarray}
  Passing to the limit and using (\ref{PW.27.07}) and that  $M_k\circ a_k=\tilde M_k\circ \tilde a_k$, we obtain
  $$  
   \gamma \|F(f)-F(g)\|_Y \geq \|f-g\|_X,
   $$
   and the proof is complete.
  \hfill $\Box$
  \nl

  The above argument brings up the interesting question of when   
  $\delta^*_{n,\gamma}(K)_X$ is  attained.     To prove a result in this direction, we recall the well-known  vector version of the Banach limit,  see Appendix C in \cite{BL}. 
   Given   $X$, let $\ell_\infty(\N,X)$  denote the space
of all  sequences $\vec f=(f_k)_{k=1}^\infty$ with $f_k\in X$ ,  $k\ge 1$, equipped with the norm
$\|\vec f\|_\infty=\sup_k\|f_k\|_X$.    The following theorem holds.
 \begin{theorem}
 \label{T:BP}
   For each Banach space $X$, there exists a
norm one linear operator
 $$
 L:\ell_\infty(\N,X)\rightarrow X^{**},
 $$
 such
that  $L(\vec f)=g$  whenever  $\vec f=(f_k)_{k=1}^\infty$, $f_k\in X$,  and $\lim_{k\to \infty} f_k=g\in X$. Note  that we have
$$
\|L(\vec f)\|_{X^{**}}\leq \limsup_{k\to \infty} \|f_k\|_X.
$$
   
  \end{theorem}

   \begin{theorem} 
      \label{T:deltaassumed}
 Let  $X$ be a separable Banach space such that  there exists  a linear norm one projection $P$  from $X^{**}$ onto  $X$. Then, 
 for every $n$ and every compact set $K\subset X$ there is a norm $\|\cdot\|_Y$ on $\R^n$ and mappings $\tilde a:K\to (\R^n,\|\cdot\|_Y)$ and
$\tilde M:(\R^n,\|\cdot\|_Y)\to X$ such that 
$$
\sup_{f\in K}\|f-\tilde M\circ \tilde a(f)\|_X=\delta_{n,\gamma}^*(K)_X.
$$
This is also the case  for   $\bar \delta_{n,\gamma}(K)_X$. 
        \end{theorem}

\noindent
{\bf Proof:}   
For each $k\ge 1$,  consider the  $\gamma$ Lipschitz maps
$a_k:K\to (\R^n,\|.\|_{Y_k})$, $M_k:(\R^n,\|.\|_{Y_k})\to X$ and the norms $\|.\|_{Y_k}$ on $\R^n$, such that
  \beqn 
\nonumber
  \lim_{k\to \infty} \sup_{f\in K}\|f-M_k(a_k(f))\|_X =\delta^*_{n,\gamma}(K)_X.
  \eeqn
We proceed  as in the proof of Theorem \ref{T:deltazero} to generate a norm $\|\cdot\|_Y$ on $\R^n$ and a sequence of mappings 
$\tilde a_k$ that converges pointwise on $K$  to the  $\gamma$ Lipschitz mapping $\tilde  a:K\to (\R^n, \|\cdot\|_Y)$ (denoted by $F$ in Theorem \ref{T:deltazero}), and a sequence of $\gamma$ Lipschitz mappings $\tilde M_k:(\R^n,\|\cdot\|_{Y_k})\to X$. Note  that 
$M_k\circ a_k=\tilde M_k\circ \tilde a_k$, and therefore
$$
\lim_{k\to \infty}\sup _{f\in K}\|f-\tilde M_k\circ \tilde a_k(f)\|_X=\delta^*_{n,\gamma}(K)_X.
$$
For $x\in  \R^n$, we consider the sequence $\overrightarrow{M(x)}:=( \tilde M_k(x))_{k=1}^\infty\in
\ell_\infty(\N,X)$  and define the mapping 
$$
M_\infty:(\R^n,\|.\|_Y)\rightarrow X^{**}
$$
 as $ M_\infty(x)=L(\overrightarrow{M(x)})$.  One easily verifies that this is  $\gamma$ Lipschitz map 
 since
\begin{eqnarray}
\nonumber
\|M_\infty(x)-M_\infty(y)\|_{X^{**}} &=&\|L(\overrightarrow{M(x)}-\overrightarrow{M(y)})\|_{X^{**}}\leq 
 \limsup_{k\to \infty} \|\tilde M_k(x)-\tilde M_k(y)\|_X\\ \nonumber
 &\leq& \gamma \limsup_{k\to \infty}\|x-y\|_k=\gamma \|x-y\|_Y.
\end{eqnarray}
 Then,    the mapping  $\tilde M:=P\circ M_\infty$,
  $ \tilde M:(\R^n,\|.\|_Y)\rightarrow X$ is a $\gamma$ Lipschitz map  since $P$ is linear projection on $X$ of norm one. For $f\in K$, we define  
  $\vec f:=(f,f,\dots)\in \ell_\infty( \N,X)$, and then
 \begin{eqnarray*}
 \|f-\tilde M\circ \tilde a(f)\|_X&=&\| P\circ  L(\vec f)-P\circ M_\infty\circ \tilde a(f)\|_X\leq 
\| L(\vec f)-M_\infty(\tilde a(f))\|_{X^{**}}\\ \nonumber
&=&\|L(\vec f)-L\big(\overrightarrow{M(\tilde a(f))}\big)\|_{X^{**}}
 \leq  \limsup_{k\to \infty}\|f-{\tilde M}_k(\tilde a(f))\|_{X}\\
 &\leq & \limsup_{k\to
\infty}      \|f-\tilde M_k(\tilde a_k(f))\|_{X}+\|{\tilde M_k(\tilde a_k(f))}-{\tilde M}_k(\tilde a(f))\|_{ X}
\\ \nonumber
&\leq & \limsup_{k\to
\infty}      \|f-{ \tilde M_k(\tilde a_k(f))}\|_{ X}+\gamma\| \tilde a_k(f)-\tilde a(f)\|_k\\
\nonumber
&\leq & \limsup_{k\to
\infty}      \|f-{ \tilde M_k(\tilde a_k(f))}\|_{X}+\gamma\|{\tilde a}_k(f)-\tilde  a(f)\|_{\ell_2^n}\\
\nonumber
  &\leq& \delta_{n,\gamma}^*(K)_X.
 \end{eqnarray*}
 Thus, we get 
 $$
 \sup_{f\in K}\|f-\tilde M\circ \tilde a(f)\|_X\leq \delta_{n,\gamma}^*(K)_X,
 $$
and the proof is completed. To show the theorem for $\bar \delta_{n,\gamma}(K)_X$,    it suffices to repeat those arguments assuming that 
$a_k$'s are defined  on $X$.  \hfill $\Box$

\begin{remark} 
Clearly a reflexive Banach space $X$  is complemented in
$X^{**}$ by a linear projection of norm one, namely the identity. The same holds for  $L_1([0,1])$ and
$L_\infty([0,1])$. However $\cC([0,1])$ is not complemented in
$\cC([0,1])^{**}$.
\end{remark}

 \section{Stable nonlinear widths   bound entropy: Carl type inequalities}
 \label{sec:bounds}

 In this section,  we  study whether  $\delta_{n,\gamma}^*(K)_X$ can  go to zero faster than the entropy numbers of $K$.  
 To understand this question, we shall prove bounds for the entropy numbers  $\e_n(K)_X$ in terms of $\delta_{n,\gamma}^*(K)_X$.   
 The inequalities we obtain  are analogous to the bounds on entropy in terms of Kolmogorov widths as given in Carl's inequality.    
 Before formulating our main theorem,  let us note that we  cannot expect  inequalities  of the form
 \be
 \label{cantbe}
 \e_n(K)_X\le C\delta^*_{\alpha n,\gamma}(K)_X,\quad n\ge 1,
 \ee
 with $\alpha>0$ a fixed constant.
 For example,   take $X=\ell_p(\N)$ with $1\leq p<\infty$ and define 
  $$
 K:= K_m:=\{(x_1,\dots,x_m,0,\dots) \quad : \quad \sum_{j=1}^m|x_j|^p\leq 1\}\subset \ell_p(\N).
  $$
  Then $\bar \delta_{n,1}(K_m)=\delta^*_{n,1}(K_m)=0$,  provided $n\ge m$.  Indeed, in this case, we can take, 
  $$
  a_n:X\to{\cnew\ell_p^n},\quad M_n: {\cnew\ell_p^n}\to X,
  $$
  where
  $a_n(x)=(x_1,\dots,x_n)$ when $x=(x_1,x_2,\dots)\in X$ and $M_n((x_1,\dots,x_n)):= (x_1,\dots,x_n,0,0,\dots)$.
  Now, given any $\alpha>0$, we choose $n$ so that $\alpha n\ge m$ and find that the right side of \eref{cantbe} is zero but the left side is not.

 \subsection{A weak inequality for entropy}   While direct inequalities like \eref{cantbe} do not hold,  
 we shall prove a weak inequality between the entropy numbers $\e_n(K)_X$ and the stable widths
 $\delta_{n,\gamma}^*(K)_X$.  To formulate our results, we  assume that $\delta_{n,\gamma}^*(K)_X\rightarrow 0$ 
 as $n\rightarrow 0$, and consider the function
 \be
 \label{defphi}
 \phi(\e):=\phi_{K,\gamma}(\e):=\min\{m:\delta_{m,\gamma}^*(K)_X\le \e\}.
 \ee
 We shall use the following lemma.

\begin{lemma}
\label{balllemma}  Let   $\gamma>0$ and let $K\subset X$ be a compact subset of the Banach space $X$.
  Let us fix a point  $f_0\in K$, $\delta>0$,  and consider the ball $B:=B(f_0,\delta)$ in $X$ of radius $\delta$, centered at $f_0$ and    $B_K:=\  K\cap B$.
 Then $B_K$ can be covered by $N$  balls of radius $\delta/2$, where
\be
\label{bl}
N\le A^m,\quad \hbox{with}\,\, A:=1+16\gamma^2,\  m:=\phi(\delta/8),
\ee
 
\end{lemma}

\noindent
{\bf Proof:} Let  $N$  be the largest number  such that there exist  points  $f_1,\dots,f_N$ from $B_K$ such that 
\be
\label{maximal}
\|f_i-f_j\|_X \ge \delta/2, \quad i\neq j.
\ee
Since  $f_1,\dots,f_N$ is a maximal number of points from $B_K$ satisfying the separation condition \eref{maximal}, it follows that any $f\in B_K$ must be in one of the balls centered at $f_j$ of radius $\delta/2 $.   
   So, we want to  bound $N$.  
 
  Let $m=\phi(\delta/8)$, namely  $m$ is  the minimal index for which  $\delta^*_{m,\gamma}(K)_X\le \delta/8$. In what follows, we assume that 
there are mappings $a_m,M_m$ for which  $\delta^*_{m,\gamma}(K)_X$ is assumed, that is 
$\|f-M_m(a_m(f))\|_X\le \delta/8$, where $a_m: K\to (\R^m,\|\cdot\|_{Y_m})$ and $M_m:(\R^m,\|\cdot\|_{Y_m})\to X$ 
are  both $ \gamma$ Lipschitz.    A similar proof, based on  limiting arguments,  holds in the case when the infimum 
$\delta^*_{m,\gamma}(K)_X$ is not assumed.

 Let us denote by $y_0:=a_m(f_0)$,  $y_j:=a_m(f_j)\in\R^m$ and $g_j:= M_m(y_j)\in X$ for $j=1,2,\dots,N$. Then, we know that
\be
\nonumber
\|f_j-g_j\|_X \le \delta^*_{m,\gamma}(K)_X\le \delta/8,\quad j=1,\dots,N,
\ee
and therefore
\be
\label{lowerbound0}
\|g_i-g_j\|_X\ge\|f_i-f_j\|_X -\|f_i-g_i\|_X -\|f_j-g_j\|_X\ge \delta/4,\quad i\neq j.
\ee
From the assumption that $M_m$ is $\gamma$ Lipschitz we have 
$$
\|g_i-g_j\|_X=\|M(y_i)-M(y_j)\|_{Y_m}\leq \gamma\|y_i-y_j\|_{Y_m}, 
$$
and therefore it follows from \eref{lowerbound0} that
\be
\label{lowerbound01}
\|y_i-y_j\|_{Y_m}\ge \frac{\delta}{4\gamma}, \quad i\neq j.
\ee
Since 
$$
\|y_0-y_j\|_{Y_m}=\|a_m(f_0)-a_m(f_j)\|_{Y_m}\leq \gamma\|f_0-f_j\|_X\leq \gamma\delta,\quad j=1,\ldots,N,
$$
all $y_j$'s, $j=1,2,\dots,N$, are in a ball  $B_Y :=B_Y(y_0,\gamma\delta )$ of radius $\gamma \delta $ and center $y_0$ in $\R^m$ with respect to the norm $\|\cdot\|_{Y_m}$.  We recall that for any $\eta>0$,  the unit ball in an $m$  dimensional Banach space can be covered by
 $(1+2/\eta)^m$
open balls of radius $\eta$, see \cite{P}, p.  63.   Therefore,    $B_Y$ can be covered by $(1+2/\eta)^m$ balls of radius 
$\eta \gamma  \delta$.  We take $\eta:=  8^{-1}\gamma ^{-2}$ so that the radius of each of these balls is $\frac{\delta }{ 8 \gamma}$.  Then, in view of \eref{lowerbound01}, each of these balls has at most one  of the points $y_j$, $j=1,\dots,N$.  This tells us that 
\be
\nonumber
N\le  (1+2/\eta)^m \le   (1+ 16 \gamma^2)^{m},
\ee
and thus proves the lemma.
\hfill $\Box$  
 \begin{theorem}
 \label{T:cgen}
 Let   $K \subset X$ be a  compact subset of the Banach space $X$  and assume $K$ is contained in a ball with radius $R$. Let  $\e>0$ and $L$ be the smallest integer such that 
 $2^L\e \ge R$. Then $K$ can be covered by $N(\e)$ balls where
 \be
 \label{cover}
N(\e) \le A^{\sum_{k=1}^{L} \phi(2^k\e/8)},\quad A:= 1+16 \gamma^2.
\ee
\end{theorem}
\noindent
 {\bf Proof:} Let $\e_k:=2^k\e$,  $k=0,1,\dots, L$, and  $m_k:=\phi(\e_k/8)$.  We know that $K$ is contained in the ball $B$ of radius $\e_L$ which without loss of generality we can assume is centered at $0$.  From Lemma \ref{balllemma}, we have that $K$ is contained in $A^{m_{L}}$ balls of  radius $\e_{L-1}$.  We can apply Lemma \ref{balllemma} to each of these new balls and find that $K$ is contained in 
 \be
\nonumber
 A^{m_{L}}\cdot A^{m_{L-1}} = A^{m_{L}+m_{L-1}}
 \ee
  balls of radius $\e_{L-2}$.  Continuing in this way, we have that $K$ is contained in $N(\e)$ balls of radius $\e=\e_0$, where
  \be
  \nonumber
  N(\e)\le A^{\sum_{k=1}^{L} m_k}.
  \ee
 This proves the theorem.  \hfill $\Box$

 \subsection{Carl type inequalities for a general Banach space}
 \label{SS: Carl}

 We can apply the last theorem to derive  bounds on entropy numbers from an assumed decay of $\delta_{n,\gamma}^*(K)_X$ in  the following way.   
From the assumed decay, we obtain bounds on the growth of   $\phi(\e)$ as $\e\to 0$ . Then we
use  these bounds in Theorem \ref{T:cgen} to derive  a bound on the number $N(\e)$ of balls of radius $\e$ needed to cover $K$. The latter   
  then translates into bounds on $\e_n(K)_X$.  We illustrate this approach with two examples
 in this section.  The first is the usual form of Carl's inequality as stated in the literature.

 \begin{theorem}
 \label{mapentthm}
 Let   $r,\gamma \geq 1$.   If $K  $ is  any compact subset of a Banach space $X$,  then  
 \be
\nonumber
\e_n(K)_X\le C (n+1)^{-r} {\displaystyle  \sup_{m\ge 0}} (m+1)^r\delta_{m,\gamma}^*(K)_X,\quad n\ge 0,
 \ee           
 with $C$ depending only on $r$ and $\gamma$.
   \end{theorem}
   
 \noindent
 {\bf Proof:}.  We fix $r>0$ and $\gamma>0$ and let 
 $$
 \Lambda:=\sup_{m\ge 0}(m+1)^r\delta_{m,\gamma}^*(K)_X.
 $$  
 If $\Lambda=\infty$, there is nothing to prove and so we assume $\Lambda <\infty$.   We claim that
    \be
  \label{boundphi}
  \phi( \Lambda 2^{- \alpha r})\le 2^\alpha,\quad \alpha\in\R.
  \ee
  Indeed,   this follows from the definition of $\phi$ and  the fact that 
 $$
 \delta_{n,\gamma}^*(K)_X\le \Lambda (n+1)^{-r},\quad n\ge 0.
 $$

  Since  $K$ is compact, it is   contained in a ball of some  radius $R$.  We now  define
  $\e:= 8 \Lambda 2^{-nr}$ 
  and let  $L$ be the smallest integer for which 
  $$
  2^L\e= \Lambda 2^{3+L-nr}\geq R,
  $$
 and  apply   Theorem \ref{T:cgen}. From \eref{boundphi}, we have
 $$
 \sum_{k=1}^{L} \phi(2^k\e/8)=\sum_{k=1}^L\phi(\Lambda 2^{k-nr})\leq
 \sum_{k=1}^{L} 2^{n-\frac{k}{r}}\leq 2^n\sum_{k=0}^{\infty} 2^{-\frac{k}{r}}=2^n(1-2^{-1/r})^{-1}.
  $$
Therefore,  it follows from \eref{cover} that
$$
 \displaystyle{N(\e)\le  A^{\sum_{k=1}^L\phi(\Lambda 2^{k-nr})}
 \leq A^{2^n(1-2^{-1/r})^{-1}}\le 2^{2^{n+c}}},
$$
 with $c$ an integer depending only on   $r$ and $\gamma$. 
  It  follows that
 $$
 \e_{2^{n+c} }(K)_X\le  8\Lambda  2^{-r n} =2^{cr}8\Lambda 2^{-(n+c)r}, \quad n\ge 0.
$$
  This proves the desired   inequality  for integers of the form $2^{n+c}$.  This can then  be extended to  
  all integers      by using  the monotonicity of $\e_n(K)_X$.
 \hfill $\Box$
 \nl

This same idea can be used to derive entropy bounds under other decay rate assumptions on 
$\delta_{n,\gamma}^*(K)_X$.  We mention just one other example to illustrate this point.  Suppose that
$$
\delta_{n,\gamma}^*(K)_X\le  \Lambda (\log_2(n+1))^\beta   (n+1)^{-r},\quad n\ge 0,
$$
for some $r>0$ and some $\beta\in\R$.
Then the above argument gives
$$
\e_n(K)_X\le C\Lambda (\log_2(n+1))^\beta  (n+1)^{-r},\quad n\ge 1,
$$
with now $C$ depending only on $r,\beta,\gamma$.

\begin{remark}
The same results obviously hold for $\o \delta_{n,\gamma}(K)_X$
since it is larger than $\delta_{n,\gamma}^*(K)_X$.
\end{remark}

 
\subsection{Carl's inequality does not hold for manifold widths}
\label{SS:Carlno}
 
It is easy to see   that  Carl's inequality does not  hold for  the manifold widths $ \delta_n(K)_X$, 
  where the assumption on the mappings $a,M$ are only that these maps are continuous.  
   For a simple example, let $X=\ell_2(\N)$ and let $(\alpha_j)_{j\ge 1}$ be any  strictly decreasing  sequence of positive numbers which tend  to $0$.   We consider the set 
   $$
    K=K(\alpha_1, \alpha_2, \ldots):=\{\alpha_je_j\}_{j\ge 1}\cup \{0\} \subset X,
   $$
  where $e_j$, $j=1,2,\dots$, is the canonical  basis for $\ell_2(\N)$. 
   {\cnew  For each $k\geq 1$ we define  
   continuous maps   $a_k:K\to  \R$ by 
$$
a_k(0)=\alpha_k, \quad a_k (\alpha_je_j)= \alpha_{\min (j,k)}, \quad j=1,2,\dots,
$$
 and $M_k:\R\to X$  as  the piecewise linear function with breakpoints $0,\alpha_k,\dots,\alpha_1$, defined by the following conditions
\begin{eqnarray*}
M_k(t)=\begin{cases}
0, \quad t\leq 0,\\
\alpha_1e_1, \quad t\geq \alpha_1,\\
\alpha_j e_j, \quad \mbox{ for } t=\alpha_j\quad j=1,\ldots,k.
\end{cases}
\end{eqnarray*}
Clearly $M_k(a_k(x))=x$ when $x=\alpha_je_j$  with  $j\le k$. For any other $x\in K$ we have $M_k(a_k(x))=\alpha_ke_k$, and so    
$$
\sup_{x\in K}\|x-M_k(a_k(x))\|_{\ell_2(\N)}=\sup_{j>k}\|\alpha_je_j-\alpha_ke_k\|_{\ell_2(\N)}<\sqrt{2} \alpha_k.
$$ 
Since $\alpha_k\rightarrow 0$ as $k\rightarrow \infty$,  we get $\delta_1(K)_{\ell_2(\N)}=0$, and thus
 $$\delta_n(K)_{\ell_2(\N)}=0 \quad  \mbox{ for }n=1,2,\dots.$$
 
Next, we bound the entropy numbers of $K$ from below.  For $ 1\le j \leq 2^n$ and any $k\neq j$, we have
$$
\|\alpha_je_j-\alpha_ke_k\|_{\ell_2(\N)}=\sqrt{\alpha_j^2+\alpha_k^2}> \alpha_{j}\geq\alpha_{2^n}.
$$
  So if we take $\e:=\frac{1}{2}\alpha_{2^n}$ with $n\ge 1$, then any attempt to cover $K$ with  $2^n$ balls with radius $\e_0\leq\e$ will fail since
  every ball in this set will contain exactly one of
  the $\alpha_je_j$, $j=1,\ldots,2^n$ and no more elements from $K$.   This gives that 
$
 \e_n(K)_{\ell_2(\N)}\ge \frac{1}{2}\alpha_{2^n}.  
$ 

We can now show that Carl's inequality cannot hold for any $r>0$.  Given such an $r$, we take for $K=K(\alpha)$ the set corresponding to a 
 sequence ${\bf \alpha}=(\alpha_1,\alpha_2,\ldots)$, where
$$
\alpha_n:=\frac{1}{[1+\log_2 n]^{r/2}}.
$$
We have that 
$$
\e_n(K(\alpha))_{\ell_2(\N)}\ge  \frac{1}{2(n+1)^{r/2}}, \quad \hbox{while}\,\,\delta_n(K(\alpha))_{\ell_2(\N)}=0\le n^{-r}, \quad n\ge 1.
$$

Finally, let} us observe that in the above construction  of $a_k,M_k$ for $K$,  the mapping $a_k$ is $1$ Lipschitz.  On the other hand, $M_k$ has poor Lipschitz constant.
Note that  since 
$$
\|M_k(\alpha_k)-M_k(\alpha_{k-1})\|_X\ge \alpha_{k-1},
$$   
 the Lipschitz constant of $M_n$  is  at least of size $ \frac{ \alpha_{k-1} }{\alpha_{k-1}-\alpha_k}$.  When $(\alpha_j)$ tends to zero slowly as in our example, then these Lipschitz constants tend to infinity.


    \subsection {Finer results on Carl's inequality}
   \label{SS:finer}
 
 Our motivation for the results in this section is the following.  One may argue that requiring   that the maps $a,M$ are Lipschitz
 is too severe and  perhaps  stability can be gained   under   weaker assumptions son these mappings.  The results of this section show
 that this is indeed the case.  Namely, we show that to establish a form of numerical stability, it is enough to have  the mapping $a$ bounded
 and the mapping $M$ satisfy  a  considerably weaker mapping property than the requirement that it be Lipschitz.  We then go on to show that
 even under these weaker assumptions on the mappings $a,M$, one can compare the error of approximation on a model class $K$ with the entropy numbers of $K$.

 Let $K$ be a compact set in the Banach space $X$ and  recall the notation
 $A:=M\circ a$ and 
 $$
 E_A(K)_X:=\sup_{f\in K} \|f-A(f)\|_X.
 $$
 We introduce  the following new properties on the pair $(a,M)$ of mappings: \vskip .1in
 \noindent
 (i)  $a: K\to (R^n,\|\cdot\|_Y)$ is bounded, i.e., $\|a(f)\|_{Y}\le\gamma \|f\|_X,\quad f\in K$;
 \vskip .1in
 \noindent
 (ii)  $M: (R^n,\|\cdot\|_Y)\to X$  satisfies 
 \be
 \label{satisfiesii}
  \|M(x)-M(y)\|_X\le  \gamma \|x-y\|_{Y}^\beta+  E_A(K)_X, \quad x,y\in\R^n,
 \ee
   where   $\gamma,\beta >0$  are fixed and $\|\cdot\|_{Y}$ is a norm  on $\R^n$.

\vskip .1in
\noindent 
Obviously,     the assumption (i) is much weaker than the assumption that $a$ is Lipschitz. Notice that (ii) is only requiring that $M$ is a Lip $\beta$ mapping for $x,y$ sufficiently far apart which is weaker that Lipschitz when $\beta\geq 1$ and stronger when $0<\beta\leq 1$.

Using these properties, we define the following {\em  bounded stable  manifold} width
\be
\label{smw}
\tilde \delta_{n,\gamma,\beta}(K)_X:= \inf_{a,M,\|\cdot\|_Y} \sup_{f\in K}\|f-M(a(f))\|_X,
\ee
where the infimum is over all maps $a,M$ satisfying (i) and (ii) and all norms $\|\cdot\|_Y$ on $\R^n$.  Clearly, we have $\delta_{n.\gamma,\beta}(K)_X\le
\delta_{n,\gamma}^*(K)_X$ for all $n\ge 1$.
We show that  properties (i) and (ii)   still guarantee a form of numerical  stability.
  \begin{theorem}
  \label{th:stability} 
  If the pair $(a,M)$ satisfies  {\rm (i)} and {\rm (ii)} with respect to the norm $\|\cdot\|_Y$ on $\R^n$ for some 
  $\beta>0$, then the approximation operator  $A:=M\circ a$ is stable in the following sense.   If in place of 
  $f\in K$ we input  $g\in K$ with $\|f-g\|_X\le \eta$ and in place of $y=a(g)$ we compute $y'$ with
  $\|y-y'\|_{Y}\leq \eta$, then 
  \be
  \label{result1}
  \|f-M(y')\|_X\le  2E_A(K)_X+  \eta+ \gamma  \eta^\beta.
  \ee
  \end{theorem}
  
  \noindent{\bf Proof:}  
  Since $A(g)=M(a(g))=M(y)$, we have  
 $$
 \| f-M(y^\prime)\|_X\leq \|f-g\|_X+\|g-A(g)\|_X+\|M(y)-M(y^\prime)\|_X\leq \eta+E_A(K)_X+E_A(K)_X+\gamma\eta^\beta,
  $$
 where we have used  (ii).\hfill $\Box$
 \nl
 
 The above theorem shows that we can obtain a form of numerical stability under rather weak assumptions on $a,M$.
 The question now is whether it is still true that when using such mappings,   the approximation error cannot go to zero
 faster than entropy numbers.  That is, do we still have a form of Carl inequality.  The following theorem shows that this is indeed the case,  up to a logarithmic loss. 
In formulating the  theorem,  we let
%
$C_0(K):=\sup_{f\in K}\|f\|_X,$
which  is finite because by assumption $K$ is compact.

 \begin{theorem}
 \label{mapentthm1}
 Let  $r,\gamma,\beta >0$.   If $K  $ is  any compact subset   of a Banach space $X$, then  

 \be
 \label{T11}
\e_{cn\ln n}(K)_X\le  (n+1)^{-r} {\displaystyle  \sup_{m\ge 0}} (m+1)^r \tilde \delta_{m,\gamma,\beta}(K)_X,\quad   n\ge 3,
 \ee           
 with $c$ depending only on $r,\beta,\gamma$ and $C_0(K)$.
   \end{theorem}

\noindent
{\bf Proof:}    Let $\tilde \delta_n:=\tilde \delta_{n,\gamma,\beta}(K)_X$, $n\ge 1$.     
We assume that the right side of \eref{T11} is finite since otherwise there is nothing to prove.  
Given any $\e>0$, we let $m=m(\e)$ be the smallest integer such that 
\be
\label{condn}
\tilde \delta_m\le   \e/4.  
\ee
  We fix for now such  a pair $(\e,m)$.  Suppose that
  $\{f_1,\dots,f_N\}$ is the largest collection  of points in $K$ such that $\|f_i-f_j\|_X\ge \e$ for all $i,j$.  Then, the 
  balls centered at the $f_j$ with radius $\e$ cover $K$.  We want now to bound $N$.
  
 Let the pair $(a_m,M_m)$ satisfies (i-ii) with respect to the norm $\|\cdot\|_{Y_m}$ and achieves the accuracy $\tilde\delta_m$ (in case the accuracy is not actually attained,
  a slight modification of the  argument below gives the result).  It follows  from \eref{condn} that the mapping $A:=A_m=M_m\circ a_m$
  satisfies 
 \be
 \label{relation}
 \tilde \delta_m=E_{A}(K)_X\le \e/4.
 \ee
 Now, consider 
 $$y_j:= a_m(f_j)\in\R^m, \quad g_j:=M_m(y_j)\in X, \quad j=1,\dots,N.
 $$   
 Because of (i),     the points $y_j$, $j=1,\dots, N$,  are all in the  ball  $B$ centered at $0$ of 
 radius   $R:=\gamma C_0(K)$ with respect to the norm $\|\cdot\|_{Y_m}$.     Since
 $$
 \|f_j-g_j\|_X=\|f_j-A(f_j))\|_X\leq \tilde \delta_m\leq \e/4,
 $$
we have that whenever $i\neq j$, 
 \be
 \label{pp}
\|M(y_i)-M(y_j)\|_X= \|g_i-g_j\|_X\ge \|f_i-f_j\|_X-\|f_i-g_i\|_X-\|f_j-g_j\|_X\geq \e/2.
\ee
Combining  condition (ii),   \eref{relation}, and \eref{pp},  we have that $y_j\in\R^m$, $j=1,\ldots,N$,  satisfy
   
  \begin{eqnarray}
\nonumber
\gamma  \|y_i-y_j\|^\beta _{Y_m}&\ge&  \|M(y_i)-M(y_j)\|_X-E_A(K)_X\geq \e/2-E_A(K)_X\\ \nonumber
&=&\e/4+(\e/4-E_A(K)_X)\geq \e/4, \quad i\neq j.
\end{eqnarray}
In other words, $\{y_1,\ldots,y_N\}$  are in the ball $B$ and they are separated in the sense that
    \be
  \label{separation1}
  \|y_i-y_j\|_{Y_m}\ge \[\frac{  \e}{4\gamma}\]^{\frac{1}{\beta}} =:\tau, \quad i\neq j.
  \ee
   We take a minimum covering of  the ball $B$  by balls  $B_1,\dots,B_M$ of radius 
$\tau/2$.
    Then, in view of \eref{separation1}, each of these balls has at most one of the points $y_j$, $j=1,\dots,N$, and therefore $N\le M$.  
    As  we have used  earlier, for any $\eta>0$,  the   unit ball  with respect to $\|\cdot\|_{Y_m}$   can be covered by
 $(1+2/\eta)^m$ balls of radius $\eta$. This tells us that 
\be
\label{boundN1}
N\le M\le     \[ \frac{C_1}{\e}\]^{m/\beta},
\ee
 with $C_1$  depending only on $\gamma, \beta$ and $C_0(K)$.
  
We can now finish the proof of the theorem.  If     $C_2:=\sup_{m\ge 0}(m+1)^r \tilde \delta_m$ is finite, we take $\e=C_2(n+1)^{-r}$. 
   We can find $c_r\in \N$, depending on $r$, such that 
$4^{1/r}(n+1)\leq c_rn+1$, for $n\geq 3$, and thus
 $$
 \tilde \delta_{c_rn}\leq \frac{C_2}{(c_rn+1)^r}\leq \frac{C_2}{4(n+1)^r}=\e/4.
 $$  
Because of the definition  $m=m(\e)$, we have that 
$m(\e)\le c_r n$.   Hence,
it follows from \eref{boundN1}  that $K$ can be
  covered with at most
  \be
  \label{boundN2}
\[ \frac{ C_1(n+1)^r}{C_2}\]^{c_r n/\beta}%
 \le 2^{cn\ln n}
  \ee
 balls of radius $\e$.  Here   $c$ depends only on   $\beta,\gamma,r$,  and $C_0(K)$.
 In other words
 $$ 
 \e_{cn\ln n}(K)_X\le C_2(n+1)^{-r},
 $$
 which is the desired result. \hfill $\Box$


 \section {Bounds for stable manifold  widths in a Hilbert space}
 \label{sec:uboundsH}
 
 The previous section gave lower bounds in terms of entropy numbers   for the optimal possible performance 
 when using Lipschitz stable approximation.  We now turn to the question of whether these performance bounds can actually  be met.
 In this section, we consider the case when the performance error is measured  in a Hilbert space  H.  
  The following theorem proves that in this case there always exits Lipschitz stable numerical algorithms 
  whose error behaves like the entropy numbers. Hence,  this result combined  with the Carl type inequalities 
shows that stable  manifold widths and entropy numbers behave the same in the case of Hilbert     spaces.

 \begin{theorem}
\label{TH2}
Let $H$ be a Hilbert space and $K\subset H$ be any compact subset of $H$.  Then for $\gamma=2$, any $n\ge 1$,  we have
\be
\label{compare}
\delta_{26n}^*(K)_H :=\delta_{ 26n,\gamma}^*(K)_H\le \bar \delta_{26n,\gamma}(K)_H\le3\e_n(K)_H.
\ee
\end{theorem}
\noindent
{\bf Proof:}
Let us fix $n$ and consider the discrete set  
$$
\cK:=\cK_n:=\{f_1,\ldots,f_{2^n}\}\subset K
$$ 
with the property  that every $f\in K$ can be approximated by an element from $\cK_n$ with accuracy 
$\varepsilon_n(K)_H$.   That is, 
for every $f\in K$ there is $f_j\in \cK$, such that 
\be
\|f-f_j\|_H\leq \varepsilon_n(K)_H.
\ee
For  the set of $2^n$ points $\cK_n\subset H$ we apply the Johnson-Lindenstrauss Lemma,  see Theorem 2.1 in  \cite{DG} for the version we use.
According to this theorem, for any $0<\e<1$, we can find a linear map   $a_\e:\cK_n\to \ell_2^{c(\e)n}$ such that
\be
\nonumber
\sqrt{\frac{1-\e}{1+\e}}\|f_i-f_j\|_H\leq \|a_\e(f_i)-a_\e(f_j)\|_{\ell_2}\leq  \|f_i-f_j\|_H, \quad i,j=1,\ldots,2^n,
\ee
whenever  $c(\e)$ is a positive integer  satisfying $c(\e)\geq \frac{4\ln 2}{\e^2/2-\e^3/3}$.

We take $\e=3/5$ and  find we can take $c(\e)=26$.  This gives a   linear map
$$
a:\cK_n\to \ell_2^{ 26n},
$$
 for which
\be
\label{JL1}
\frac{1}{2}\|f_i-f_j\|_H\leq \|a(f_i)-a(f_j)\|_{\ell_2}\leq  \|f_i-f_j\|_H, \quad i,j=1,\ldots,2^n.
\ee
Using the Kirszbraun extension theorem,  see Theorem 1.12
from \cite{BL}, page 18, the mapping $a$ can be extended from $\cK_n$ to the whole $H$  preserving the Lipschitz constant $1$. 
Let us denote by $\cM_n$ the image  of $\cK_n$ under $a$, that is the discrete set 
$$
\cM_n:=\{a(f_j):\,f_j\in \cK_n\}\subset \R^{26n}.
$$   
Now   consider the map $M:(\cM_n,\|\cdot\|_{\ell_2})\rightarrow H$, defined by
$$
M(a(f_j))=f_j,\quad j=1,\ldots,2^n.
$$
Clearly 
$$
\|M(a(f_i))-M(a(f_j))\|_H=\|f_i-f_j\|_H\leq 2\|a(f_i)-a(f_j)\|_{\ell_2},
$$
and therefore $M$ is a Lipschitz map with a Lipschitz constant $2$. According to the Kirszbraun extension theorem, we can 
extend $M$ to a Lipschitz map on the whole $\ell_2^{26n}$ with the same Lipschitz constant $2$.

Let us now consider the approximation algorithm   $ A$ defined by   $A:=M\circ a$.  If 
  $f\in K$, there is  an  $f_j\in \cK_n$, such that 
$\|f-f_j\|_H\leq \e_n(K)_H$.  Therefore,
\be
\nonumber
f-A(f)=(f-f_j)+(f_j-M(a(f_j)))+(M(a(f_j))-M(a(f))),
\ee
and since $f_j=M(a(f_j))$, we have that 
\begin{eqnarray}
\nonumber
\|f-A(f)\|_H&\leq &\|f-f_j\|_H+\|M(a(f_j))-M(a(f))\|_H\leq \varepsilon_n(K)_H+2\|a(f)-a(f_j)\|_{\ell_2}\\
\nonumber
&\leq& \e_n(K)_H+2\|f-f_j\|_H\leq 3\varepsilon_n(K)_H,
\end{eqnarray}
 which proves the theorem.\hfill $\Box$

 We can combine the last result with the results of the previous section to obtain the following corollary.
 
 \begin{cor}
 \label{C1}
 Let $\gamma\geq 2$. If $K\subset H$ is a compact set in a Hilbert space $H$ and if $r>0$, then 
 $$
 \delta_{n,\gamma}^*(K) _H={\cal O}((n+1)^{-r}), \quad n\ge 0,\quad \hbox{if and only if}\quad  \e_n(K)_H={\cal O}((n+1)^{-r}), \quad n\ge 0.
 $$  
 The same result holds if  $\delta_{n,\gamma}^*$ is replaced by $\bar \delta_{n,\gamma}$.
 \end{cor}


 \section{Comparisons for an arbitrary  Banach space $X$}
 \label{sec:Banach}
  
 In this section, we consider bounding the stable  manifold widths by entropy numbers in the case of an arbitrary Banach space $X$.  
 Let us note that for such a general Banach space we  can no longer have a direct bound for $\bar \delta_{n,\gamma}(K)_X$ in terms of entropy numbers.   
 Indeed, for any compact set $K$ and any Banach space, the entropy numbers of $K$ tend to zero.   However, we know that $\bar \delta_{n,\gamma}(K)_X$
 tends to zero for all compact sets $K$ only if
 $X$ has the $\gamma^2$-bounded approximation property, see Theorem \ref{T:tozero}. Since there are Banach spaces without this property, we must expect a 
 loss when compared to the theorems of the
previous section.  We present in this section  results that exhibit a loss in both   the growth of the Lipschitz constants and in the rate of decay 
of   $\bar \delta_{n,\gamma}(K)_X$, as $n$ tends to infinity.
It is quite possible that the results of this section may be improved with a deeper analysis.

 \begin{theorem}
\label{BanachTH2}
Let $X$ be a Banach space and $K\subset X$ be a compact subset of $X$. Then,  
there  is a  fixed  positive constant $C$,  such that  for each $n\geq 1$  there are  Lipschitz mappings
$$
a_n:X\to (\R^{26n},\|\cdot\|_{\ell_{\infty}}), \quad M_n: (\R^{26n},\|\cdot\|_{\ell_{\infty}})\to X
$$  whose Lipschitz constants are at most $Cn^{5/4} $  and   
\be
\label{EQ1}
\sup_{f\in K}\|f-M_n(a_n(f))\|_X\leq Cn^{5/2} \varepsilon_{n}(K)_X,\quad  n=1,2,\ldots.
\ee
\end{theorem}

\noindent
{\bf Proof:}    
As in the proof of Theorem \ref{TH2}, we fix $n>0$, and consider the discrete set  
$$
\cK_n:=\{f_1,\ldots,f_{2^n}\}\subset K,
$$ 
with the property  that for every $f\in K$ there is $f_j\in \cK_n$, such that 
\be
\|f-f_j\|_X\leq \varepsilon_n(K)_X.
\ee
For the discrete set $\cK_n\subset X$ of $2^n$ points we apply Proposition 1 from \cite{B},
according to which we can construct a 
bi-Lipschitz map $\tilde a_n$ from $\cK_n$  into a Hilbert space $H$, 
$$
\tilde a_n:(\cK_n,\|\cdot\|_X)\to H,\quad \tilde a_n^{-1}:(\cH_n,\|\cdot\|_H)\to \cK_n, \quad \hbox{where}\quad \cH_n:=\tilde a_n( \cK_n)\subset H,
$$  
such that $\tilde a_n$  is $C_1 n^{5/4}$ Lipschitz map and $\tilde a_n^{-1}$ is $C_2 n^{-1/4}$ Lipschitz map.
Using the version of the  Johnson-Lindenstrauss lemma  as in the proof of Theorem \ref{TH2}, we get a map 
$$
J:(\H_n,\|\cdot\|_H)\to \ell_2^{26n}, 
$$
such that $J$ and $J^{-1}$ are $2$ Lipschitz maps.
We also consider the identity map
$$
I: \ell_2^{26n}\to \ell_\infty^{26n},
$$
where $I$ is $1$ Lipschitz map and $I^{-1}$ is $\sqrt {26n}$ Lipschitz map. Thus, the map
$$
a_n:=I\circ J\circ \tilde a_n:\cK_n\to \ell_\infty^{26n}
$$
 is a $Cn^{5/4}$ Lipschitz map which,   see \cite[Lemma 1.1]{BL}, can be extended to a map 
 $$
 a_n:X\to \ell_\infty^{26n}
 $$  
 with the same Lipschitz constant. 

Next, we proceed with the construction of $M_n$.
First, we denote by $\cM_n\subset \R^{26n}$ the image  of $\cK_n$ under $a_n$, that is the discrete set 
$$
\cM_n:=\{a_n(f_j):\,f_j\in \cK_n, \,j=1,\ldots,2^n\}\subset \R^{26n},
$$ 
and consider the map
$$
\tilde M_n:=\tilde a_n^{-1}\circ J^{-1}\circ I^{-1}:(\cM_n,\|\cdot\|_{\ell_\infty})\to X.
$$
 From the above observations it follows that $\tilde M_n$ is a $Cn^{1/4}$ Lipschitz map.
According to Theorem 1 from \cite{JLS}, we can extend $\tilde M_n$ to a Lipschitz map $M_n$ from $\ell_\infty^{26n}$ into $X$ with the Lipschitz constant $Cn^{5/4}$.

Now that $a_n$ and $M_n$ are constructed, we continue  with the analysis of the approximation power of the mapping $M_n\circ a_n$.
 We fix $f\in K$, find $f_j\in \cK_n$, such that 
$\|f-f_j\|\leq \varepsilon _n(K)_X$. 
Clearly,
\begin{eqnarray}
\nonumber
\|f-M_n\circ a_n(f)\|_X&\leq &\|f-f_j\|_X+\|M_n(a_n(f_j))-M_n(a_n(f))\|_X
\\ \nonumber
&\leq &\varepsilon_n(K)+C n^{5/4}\|a_n(f)-a_n(f_j)\|_{\ell_\infty }\\
\nonumber
&\leq& \varepsilon_n(K)+Cn^{5/2}\|f-f_j\|_X \leq Cn^{5/2}\varepsilon_n(K)_X.
\end{eqnarray}
Therefore, for the $Cn^{5/4}$ Lipschitz mappings $a_n$ and 
$M_n$, we have 
$$
\sup_{f\in K}\|f-M_n\circ a_n(f)\|_X\leq Cn^{5/2} \varepsilon_n(K)_X.
$$
This completes the proof. \hfill $\Box$

\begin{remark}
\label{remark:Bspace}
If we have  additional information about the Banach space $X$,  we can get better estimates than \eref{EQ1}, as illustrated in the next lemmas.
\end{remark}
\begin{lemma}
\label{L1}
 Let the  Banach space $X$ be  isometric to $\ell_\infty(\Gamma)$ for some set $\Gamma$. Then, there  is a  fixed  positive constant $C$,  
such that  for each $n\geq 1$ there are $Cn^{3/4}$ Lipschitz mappings 
$$
a_n:X\to \ell_\infty^{26n}, \quad  M_n:\ell_\infty^{26n}\to X,
$$
 with the property 
$$
 \sup_{f\in K}\|f-M_n(a_n(f))\|_{X}\leq Cn^{3/2} \varepsilon_{n}(K)_{X},\quad  n=1,2,\ldots.
 $$
\end{lemma}

\noindent
{\bf Proof:} For $\cK_n$  as in the proof of Theorem \ref{BanachTH2} and $H$ a Hilbert space, using  \cite{B},
we  construct mappings 
$$
\tilde a_n:\cK_n\to H,\quad \tilde a_n^{-1}:\cH_n\to \cK_n, \quad \hbox{where}\quad \cH_n:=\tilde a( \cK_n)\subset H,
$$
where $\tilde a_n$ is Lipschitz with constant $C_1n^{3/4}$ and $ \tilde a_n^{-1}$ with a Lipschitz constant $C_2n^{1/4}$. 
Then, with $I$ and $J$ are as in Theorem \ref{BanachTH2}, the mapping
$$
I\circ J\circ \tilde a_n:\cK_n\to \ell_\infty^{26n},
$$
 is a $C_1n^{3/4}$ Lipschitz.   We extend it to a mapping $a_n$ on the whole $X$ with the same Lipschitz constant. 

Next, we  consider 
$$
\tilde M_n:=\tilde a_n^{-1}\circ J^{-1}\circ I^{-1}:(\cM_n,\|\cdot\|_{\ell_\infty})\to X, 
$$
which is $C_2n^{3/4}$ Lipschitz.  Now, according to Lemma 1.1 from \cite{BL}, since $X$ is isometric to $\ell_\infty(\Gamma)$ for some $\Gamma$,
$\tilde M_n$ can be extended to 
$$
M_n:\ell_\infty^{26n}\to X
$$ 
with the same Lipschitz constant
$C_2n^{3/4}$.  
Then $A_n=M_n\circ a_n$ is $Cn^{3/2}$ Lipschitz, and 
\begin{eqnarray}
\nonumber
\|f-M_n\circ a_n(f)\|_X&\leq &\|f-f_j\|_X+\|M_n(a_n(f_j))-M_n(a_n(f))\|_X
\\ \nonumber
&\leq& \varepsilon_n(K)+Cn^{3/2}\|f-f_j\|_X \leq Cn^{3/2}\varepsilon_n(K)_X,
\end{eqnarray}
which gives 
$$
\bar\delta_{26n,Cn^{3/4}}(K)_{X}\leq Cn^{3/2} \varepsilon_n(K)_{X}.
$$
\hfill $\Box$

\begin{cor}\label{L2} Let $\cC(S)$ be the  Banach space of continuous functions on a  compact subset $S$ of a metric space.  Further,  let $K\subset \cC(S)$  be a compact set. Then, there  is a  fixed  positive constant $C$,  
such that  for each $n\geq 1$ there are $Cn^{3/4}$ Lipschitz mappings 
$$
a_n:\cC(S) \to \ell_2^{26n}, \quad  M_n:\ell_2^{26n}\to  \cC(S),
$$
 with the property 
$$
 \sup_{f\in K}\|f-M_n(a_n(f))\|_{\cC(S)}\leq Cn^{3/2} \varepsilon_{n}(K)_{\cC(S)},\quad  n=1,2,\ldots.
 $$
\end{cor}
\noindent
{\bf Proof:} Let us fix arbitrary $\e>0$.  Since $\cC(S)$  is separable, it follows from  \cite{MiPe}  that there exists a finite dimensional subspace $X\subset \cC(S)$ isometric to $\ell_\infty(\Gamma)$ and a linear projection $P:\cC(S)\to X$ from $\cC(S)$ onto $X$ of norm $1$ such that 
$$
\sup_{f\in K}\|f-P(f)\|\leq \e.
$$
We apply Lemma \ref{L1} to the space $X$ and its compact subset $P(K)$, according to which there are $Cn^{3/4}$ Lipschitz mappings 
$$
a_n:X\to \ell_\infty^{26n}, \quad  M_n:\ell_\infty^{26n}\to X,
$$
 with the property 
$$
 \sup_{g\in P(K)}\|g-M_n(a_n(g))\|_{\cC(S)}\leq Cn^{3/2} \varepsilon_{n}(P(K))_{\cC(S)},\quad  n=1,2,\ldots.
 $$
We next  define $\tilde a_n:\cC(S)\to \ell_\infty^{26n}$, and 
$\tilde M_n:\ell_\infty^{26n}\to \cC(S)$, 
where 
$$
\tilde a_n:=a_n\circ P, \quad \tilde M_n:=I\circ M_n,
$$
 with  $I:X\to \cC(S)$ the identity embedding from $X$ into $\cC(S)$. Clearly 
 $\tilde a_n$ and $\tilde M_n$  are both $Cn^{3/4}$ Lipschitz  and
\begin{eqnarray}
\sup_{f\in K}\|f-\tilde M_n(\tilde a_n(f))\|_{\cC(S)}&=&\sup_{f\in K}\|f-M_n(a_n(P(f))\|_{\cC(S)}\nonumber\\
&\leq&\sup_{f\in K}\left(\|f-P(f)\|_{\cC(S)}+\|P(f)-M_n(a_n(P(f))\|_{\cC(S)}\right) \nonumber\\
&\leq& \e+Cn^{3/2}\e_n(K)_{\cC(S)},
\nonumber
\end{eqnarray}
where we have used  that $\e_n(P(K))_{\cC(S)}\leq \e_n(K)_{\cC(S)}$. 
Since $\e$ is arbitrary we get the claim.
\hfill $\Box$

\section{Examples of linear and nonlinear approximation}
\label{S: examples}
 Next, we discuss a few  standard examples of approximation from the viewpoint of stable manifold widths.
 
 \subsection{Linear approximation}
 \label{SS:linearapprox}
 Let $X$ be a  Banach space, $K\subset X$ be compact,  and let $X_n$ be a linear subspace of $X$ of dimension $n$.  Let us consider approximation procedures $f\to A(f)=M\circ a(f)$ given by maps $a,M$, where 
$$
a:X\to \R^n, \quad M:\R^n\to X_n\subset X.
$$  
 If we are interested   only in such approximation methods given by continuous mappings then it is easy to see that  by using coverings and partitions of unity (see  Theorem 2.1 in  \cite{DHM}) one can achieve an approximation error for $K$ equivalent to the error $\dist(K,X_n)_X$.   Thus, $\delta_n(K)_X$ can be bounded by the $Cd_n(K)_X$ where $d_n$ is the Kolmogorov width.  The situation becomes more subtle when we require Lipschitz continuity of the mappings as we
  now discuss.

  Let   $\Phi:=\{\phi_1,\dots,\phi_n\}$ be  any basis for $X_n$ and let us consider the norm on $\R^n$, induced by the basis $\phi_1,\dots,\phi_n$, namely
 \be
\label{Rnnorm} 
\|y\|_{ Y}:=\|\sum_{j=1}^n y_j\phi_j\|_X, \quad y\in \R^n.
\ee
We define the mapping $M:(\R^n,\|\cdot\|_\Phi)\to X$, as
$$
M(y):=\sum_{j=1}^n y_j\phi_j\in X_n\subset X, \quad y\in \R^n.
$$
Clearly, $M$ is a linear mapping with norm one, and hence a 1 Lipschitz mapping. Thus, the main question is whether we can construct a mapping
 $a:X\to (\R^n,\|\cdot|\|_{\anew Y})$ that is Lipschitz.

If $X_n$ admits a bounded projection $P_n:X\to X_n$, then we can write for $f\in X$, 
 \be
  \label{proj}
  P_n(f)=\sum_{j=1}^n a_j(f)\phi_j,
  \ee
and therefore define $a$ as
$$
a(f)=(a_1(f),\dots,a_n(f))\in \R^n.
$$
 Since
$$
\|a(f)-a(g)\|_{ Y}=\|\sum_{j=1}^n a_j(f-g)\phi_j\|_X=\|P_n(f-g)\|_X\leq \|P_n\|\|f-g\|_X,
$$
$a$ is a  $\gamma_n$-Lipschitz mapping with $\gamma_n:=\|P_n\|\geq 1$.
We thus have 
$$
\bar \delta_{n,\gamma_n}(K)_X\leq \sup_{f\in K}\|f-M(a(f))\|_X=\sup_{f\in K}\|f-P_n(f)\|_X.
$$
If $X=H$ is a Hilbert space then we know there is always a projection with norm one and hence 
\be
\label{Hilbert1}
 \bar \delta_{n,1}(K)_H\le d_n(K)_H,  \quad n\ge 1.
 \ee
 For non-Hilbertian Banach spaces  every finite dimensional space admits a projection, however the norm may depend on $n$. The   Kadec-Snobar
 theorem  guarantees that there is a projection with norm $\sqrt{n}$ and so we obtain for a general Banach space $X$ and compact $K\subset X$ the bound
  \be
\label{Banach1}
 \bar \delta_{n,\sqrt{n}}(K)_X\le d_n(K)_X,  \quad n\ge 1.
 \ee
 Of course, we already know from our earlier results that relate the decay of $\bar \delta_{n,\gamma}(K)_X$ to the bounded approximation property that some growth factor is needed.    If we assume additional structure on $X$ then the quantitative growth can be better controlled.   For example, for $X= L_p$, $1<p< \infty$, 
 we can replace $\sqrt{n}$ in \eref{Banach1} by $n^{|1/2-1/p|}$, see e.g. \cite[III.B.10.]{W}.

\subsection{Compressed Sensing}
\label{SS:CS}
One of the primary settings where nonlinear approximation methods prevail is in compressed sensing which is concerned with the numerical recovery of sparse  signals.
The standard setting of compressed sensing is the following.  We consider vectors $x\in\R^N$ where $N$ is large.  Such a  vector $x$ is said to be $k$ sparse if at most $k$ of its coordinates are nonzero.  Let $\Sigma_k$ denote the set of all $k$ sparse vectors in $\R^N$. 
The goal of compressed sensing is to make a small number of $n$ linear measurements
of a vector $x$ which can then be used to approximate $x$.  The linear measurements take the form of inner products of $x$ with vectors $\phi_1,\dots,\phi_n$.  These measurements can be represented as the application of a compressed sensing matrix $\Phi\in \R^{n\times N}$  to $x$, where the  rows of $\Phi$ are the vectors 
$\phi_1,\dots,\phi_n$.

A fundamental assumption about the measurements used in compressed sensing is the so called {\it restricted isometry property of order $k$}, RIP($k,\delta_k$). 
We say that the matrix $\Phi$  satisfies the RIP($k,\delta_k$), $0<\delta_k <1$, if
\be \label{RIP}
(1-\delta_k)\|x\|_{\ell_2^N}\leq \|\Phi(x)\|_{\ell_2^n} \leq (1+\delta_k)\|x\|_{\ell_2^N}, \quad \hbox{for all}\quad x\in \Sigma_k.
\ee
 
A decoder is a mapping $M$ which takes  the measurement vector $y=\Phi(x)$ and maps it back into $\R^N$.  The vector $M(\Phi(x))$ is the approximation to $x$.
Thus, compressed sensing falls into our paradigm of nonlinear approximation as given by the two mapping 
$a:\R^N \to ( \R^n,\|\cdot\|_Y)$ with $a(x):=\Phi(x)$ and the mapping 
$M:(\R^n,\|\cdot\|_{Y})\to \R^N$.  Note that the mapping $a$ is rather special since it is assumed to be linear.

The first goal of compressed sensing is to find such mappings  for which  $M(a(x))=x$ whenever $x$ is in $\Sigma_k$.  It is easy to see that $n=2k$ is the smallest number of measurements for which this is true and it is easy to characterize all of the mappings
$a=\Phi$ that do the job (see e.g. \cite{CDD:CS}).  However, these matrices $\Phi$  and perfect reconstruction maps $M$ with  $n=2k$  are deemed unsatisfactory because   of their instability.
To discuss this and other issues connected with  compressed sensing using the viewpoint of this paper, we need to introduce a norm on $\R^N$ in which we shall measure performance.  We consider the 
 $\ell_p$  norms for $1\le p\le 2$ in what follows, therefore taking $X:=(\R^N,\|\cdot\|_{\ell_p})$.

 There are two flavors of results one can ask for in the context of compressed sensing or sparse recovery.  The strongest guarantees are in the form of  {\it instance optimality}.  To formulate this let $x\in \R^N$ and define
\be
\label{sigman}
\sigma_k(x)_p:=\inf_{y\in\Sigma_k}  \|x-y\|_{\ell_p}
\ee
to be its error of best approximation by $k$ sparse vectors.  We say that the measurement system $(\Phi,M)$ is $C$   {\it  instance optimal  of order} $k$ if
\be
\label{io}
\|x-M(\Phi(x))\|_{\ell_p}\le C\sigma_k(x)_p,\quad x\in \R^N.
\ee
A central issue in compressed sensing is how large must the number of measurements $n$ be to guarantee instance optimality of order $k$ with a reasonable constant $C$. 
  It is known, see \cite{CDD:CS}, 
that  for $p=1$, linear mappings 
 $\Phi$ based on $n$ measurements and satisfying the $ \RIP(3k, \delta_{3k})$, with $\delta_{3k}\leq \delta<(\sqrt{2}-1)^2/3$,  and the recovery map $M$ based on $\ell_1$ minimization
 \be
 \nonumber
 M(y):={\rm argmin}\{\|x\|_{\ell_1} \,:\, \Phi x=y\},
 \ee
provide instance optimality.  One can construct such matrices when $n\ge c k\log(N/k)$ with a suitable constant $c$ independent of $k$.    
 On the other hand,  see \cite{CDD:CS}, when $1<p\leq 2$, the number of measurements $n$ must necessarily grow as a power of $N$ in order to guarantee that the  instance optimality \eref{io} is achieved.  In particular, for $p=2$, instance optimality cannot hold unless $n$ is proportional to $N$.
 
 A weaker notion of performance is to consider only distortion on compact subsets $K$ of  $\R^N$.  The distortion is now measured in the worst error described by
 \be
 \label{worst}
 E(K,\Phi,M)_p:= \sup_{x\in K}\|x-M(\Phi(x))\|_{\ell_p}.
 \ee
A common family of model classes  are  the unit balls $K_q$,
 $$
 K_q:=\{x\in\R^N:\,\,\|x\|_{\ell_q}\le 1\}, \quad q<p.
 $$
 By utilizing the above results on instance optimality for $p=1$, one can derive estimates for the above error when using  a suitably chosen  compressed sensing matrix   $\Phi$ for encoding and with 
 $\ell_1$ minimization decoding $M$.  Given $p\ge 1$, one can derive bounds for the above error  for a certain range of  $q$ and show these are optimal
 by comparing this error with Gelfand widths.  We refer the reader to \cite{CDD:CS} for details.

 Our main goal  in this paper is not to restrict the measurement map $a$ to be linear but rather impose only that it is Lipschitz. 
 By relaxing the condition on $a$ to only be Lipschitz we will derive improved approximation error bounds.  We first observe  that the matrices $\Phi$, which are the canonical measurement maps  of compressed sensing,  have rather big   
  Lipschitz constants  when considered as mapping 
  from $\ell_p^N$ to $\ell_2^n$. Let us denote by  $\|\Phi\|_{\ell_p^N\to \ell_2^n}$  the norm of $\Phi$. Then the following lemma holds.

 \begin{lemma} 
If the  matrix  $\Phi$ satisfies the $\RIP(1,\delta)$,  
then for all $1\leq p\leq 2$, 
     \be 
     \label{RIPlipschitz}
     (1-\delta)^{-1} n^{-1/2} N^{1-1/p} \leq \|\Phi\|_{\ell_p^N\to \ell_2^n}
     \leq (1+\delta)  N^{1-1/p}.
     \ee
 \end{lemma}
 
 \noindent
 {\bf Proof:}  Let $\Phi:=(a_{i,j})\in \R^{n\times N}$. It follows from the $\RIP(1,\delta)$ that for $j=1, \ldots,N$,
\be
\label{po}
(1-\delta)^2\leq \sum_{i=1}^n |a_{i,j}|^2 \leq (1+\delta)^2,
\ee
and therefore
\be
\label{op1}
(1-\delta)\sqrt N \leq \|\Phi\|_{F}\leq (1+\delta) \sqrt N,
\ee
where $\|\Phi\|_{F}$ is the  Frobenious norm of $\Phi$.
 Since
  $$
  \frac{1}{ \sqrt n} \|\Phi\|_F \le \|\Phi\|_{\ell_2^N\to \ell_2^n} \leq \|\Phi\|_{F},$$
  it follows from \eref{op1} that 
\be
\label{ll}
  (1-\delta) n^{-1/2}\sqrt N \leq \|\Phi\|_{\ell_2^N\to \ell_2^n} \leq (1+\delta)\sqrt N.
\ee
We now derive bounds for  
$\Phi$ on the $\ell_p^N$ spaces, $1\le p<2$.   Let $e_j:=(0,\ldots,1,0, \ldots,0)\in \R^N$, be the $j$-th standard basis element.   We have $\|e_j\|_{\ell_1^N}=1$ and
 $$
\|\Phi e_j\|^2_{ \ell_2^n} =\sum_{i=1}^n |a_{i,j}|^2\leq  (1+\delta)^2,\quad j=1,\dots,N,
 $$
where we have used \eref{po}. Thus, for every $x=\sum_{j=1}^nx_je_j\in \ell_1$,
$$
\|\Phi x\|_{ \ell_2^n}\leq \sum_{j=1}^n|x_j|\|\Phi e_j\|_{ \ell_2^n}\leq (1+\delta)\|x\|_{ \ell_1^N}. 
$$
 In other words,
 $$
 \|\Phi\|_{\ell_1^N\to \ell_2^n} \leq (1+\delta),
 $$
 and from \eref{ll} and the Riesz-Thorin theorem we get the right inequality in \eref{RIPlipschitz}. 

 To prove the left   inequality  in \eref{RIPlipschitz} ,  we observe that from \eref{op1} there exists $1\leq i_0 \leq n$ such that 
 $$
 \sum_{j=1}^N a_{i_0,j}^2\geq \frac{ N}{(1-\delta)^2n}.
 $$
 We define $a^*:=(a_{i_0,1}, \ldots,a_{i_0,N})\in\R^N$ and  $x^*:=a^*/\|a^*\|_{\ell_2^N}$       Then we have 
 $$
 \frac{  N^{1/2}n^{-1/2}}{(1-\delta)}  \leq \sum_{j=1}^N x^*_j a_{i_0,j}=[\Phi x^*]_{i_0} \leq \|\Phi x^*\|_{\ell_2^n}\leq \|\Phi\|_{\ell_p^N\to \ell_2^n}\|x^*\|_{\ell_p^N}.
 $$
 Since $\|x^*\|_{\ell_p^N}\leq N^{1/p-1/2}$ we get the left inequality in\eref{RIPlipschitz}.
\hfill $\Box$
 \nl
 
 Since the mapping $\Phi$ is linear, its norm is the same as its Lipschitz constant.  So the above lemma shows that this Lipschitz constant
 is large, at least when we choose the norm on $\R^n$ to be the $\ell_2$ norm.  Choosing another norm on $\R^n$ cannot help much because of norm equivalences on $\R^n$ and changing norms will change the Lip constant for  the recovery mapping $M$.  We next want to show that dropping the requirement
 that $a$ is linear, and replacing it by requiring only that it is Lipschitz, dramatically improves matters.
 For now, we illustrate this only in one setting.  We consider instance optimality in $\ell_2$ which we recall fails to hold in the classical setting of compressed sensing.

Let  $X=\ell_2^N$ and let $\Phi$ be an $n\times N$ matrix which satisfies the RIP of order $2k$ (with suitable RIP constants).    Define 
$a:\Sigma_k\to \ell_2^n$ by
$$
a(x):=\Phi(x), \quad x\in \Sigma_k.
$$
It follows from the RIP that $\|\Phi  x\|_{\ell_2^n}\le C\|x\|_{\ell_2^N}$, for all $x\in\Sigma_{2k}$, and so $a$ is $C$ Lipschitz on $\Sigma_k$.   By the  Kirszbraun extension  theorem, $a$  has a $C$  Lipschitz  extension to all of $X$ which extension we continue to denote by $a$.    Note that $a$ will not be linear on $X$.

Now consider the construction of a recovery map $M$.  There is a  1 Lipschitz  inverse mapping
$M:a(\Sigma_k)\to X$ such that $M(a(x))=x$ when $x\in\Sigma_k$ (for example $\ell_1$ minimization provides such an $M$).  Again by  the  Kirszbraun extension theorem, $M$ has a 1 Lipschitz  extension to all of $\ell_2^n$, which we continue to denote by $M$.  

These new mappings
\be
\label{newmaps}
a:\ell_2(\R^N)\to\ell_2(\R^n),\quad M:\ell_2(\R^n)\to \ell_2(\R^N),
\ee
have Lipschitz constant at most $C$ for $a$ and one for $M$.  Moreover, when applied to any $x\in\Sigma_k$,  we still have $M(a(x))=x$.

Now, consider the performance of these mappings on all of $\ell_2^N$.  Given $x\in\R^N$, 
we can write $x=x_0+e$, where $x_0$ is a best approximation to $x$ from $\Sigma_k$ and $\|e\|_{\ell_2}=\sigma_k(x)_{\ell_2}$.  We have that
\begin{eqnarray}
\label{perf}
\|x-M(a(x))\|_{\ell_2}&\le& \|x_0+e-M(a(x_0))\|_{\ell_2} + \|M(a(x_0))-M(a(x))\|_{\ell_2} \nonumber \\
&\le&  \|e||_{\ell_2} +C\|e\|_{\ell_2}=(C+1)\sigma_k(x)_{\ell_2},
\end{eqnarray}
because $M(a(x_0))=x_0$ and because the composition mapping $M\circ a$ is  $C$ Lipschitz mapping.
 Thus, instance optimality can be achieved in $\ell_2$,
for $n$ of the order of $k$ up to logarithmic factors
provided one generalizes the notion of measurement maps to be nonlinear but Lipschitz, while linear measurements would impose that $n$ is of the order of $N$.

\subsection{Neural networks}
\label{SS:neuralnetworks}

This is now a very active area of research.  A neural network is a vehicle for creating multivariate functions which depend on a fixed number $n$ of parameters given 
by the weights and biases of the network.     We  consider all networks with $n$ parameters with perhaps some  user prescribed restrictions imposed on the architecture of the network.  Let us denote by $\Upsilon_n$ the outputs of such networks.    Thus the elements in $\Upsilon_n$ are multivariate functions, say with $d$ variables,  described by $n$ parameters and hence are a nonlinear manifold depending on $n$ parameters. 
 
Let us fix  a function  norm $\|\cdot\|_X$ to measure error.    Given a target function $f\in X$ (or  data observations of $f$ such as point values),  one determines the $n$ parameters $a(f)=(a_1(f),\dots,a_n(f))$ of the network which will be used to approximate $f$.
 These parameters determine the output function $M(a)$ from $\Upsilon_n$.  The decoder $M$ is explicit and simple to describe from the assumed architecture.
 For example, for the ReLU activation function this output is a piecewise linear function.   Thus, neural networks provide an approximation procedure $A(f):=M(a(f))$ of the type studied in this paper.
 
 There are by now several papers addressing the approximation properties of neural networks (see \cite{DDFHP} and the references therein).  In some cases, they advertise some surprising results.  We mention here only the results on approximating univariate  $1$ Lipschitz  functions with respect to  an $L_p$ norm on an interval $[0,1]$ by 
 neural networks with a ReLU activation function. It is shown in \cite{Shen} (with earlier results in \cite{Y}) that any function in the unit ball $K$ of Lip $1$
 can be approximated to accuracy $Cn^{-2}$ by elements from $\Upsilon_n$.   This result is on first glance quite surprising since the entropy   number $\e_n(K)_{L_p}\ge cn^{-1}$ with $c$ an absolute constant.

 So, how should we evaluate such a result? The first thing we should note is that if we view such a neural network approximation as simply a manifold approximation, then the result is not surprising.  Indeed, we could equally well construct  a one  parameter (space filling) manifold  (even with piecewise linear manifold elements) and  achieve arbitrary  approximation error for $K$.  Such a one parameter manifold is not very useful, since given $f$ or data for $f$, it would be essentially impossible to
 numerically find an approximant from the manifold with this error.  So the main issues center around the properties of $a$ and $M$.  If we require the rather minimal condition that $a$ and $M$ are continuous, we can
 never achieve accuracy better than $cn^{-1}$ in approximating the elements of $K$ using an $n$ parameter manifold as is proved in \cite{DKLT}.  We can even lessen the requirement that $a$ be continuous to just requiring
 that $a$ is bounded if we impose a little smoothness on $M$ (see Theorem \ref{mapentthm1}).  So, to achieve a rate of approximation better than $O(n^{-1})$ for $K$ using $n$ parameter neural networks, one must
 necessarily use mappings which are not continuous, even $a$ has to be poorly bounded (with bounds growing with $n$).  The question is  the numerical cost to find good parameters and whether the  numerical
 procedure to find these parameters is stable.   The results of the present paper clarify these issues.
 
 In practice, the parameters of the neural network are found from given data observations of $f$,  by typically using stochastic gradient descent algorithms with respect to a chosen loss function related to fitting the data.
 Unfortunately, there is no clear analysis of the convergence of these decent algorithms  for such optimization problems,
   although it seems to be recognized that one needs to impose constraints on the size of the steps in each iteration
 that tend to zero as the number of steps increase.   The results of the present paper may provide a better understanding of what conditions need to be imposed in the descent and what approximation results can be obtained under such constraints.

\subsection{Conclusion}
   A general question, which is  not  answered  in this paper
is to determine the asymptotic behavior of $\delta_{n,\gamma}^*(K)_X$ for classical model classes $K$ in classical Banach spaces $X$.  For example, we do not know the decay rate of $\delta_{n,\gamma}^*(K)_X$  for all of the  Besov or  Sobolev balls $K$ that compactly embed into   $L_p$,   $1\le p\le\infty$.   The asymptotic decay of these widths remains an open fundamental question. In the case that this ball is a compact subset of $L_p$, then it  is known, 
 see Theorem 1.1 in \cite{CDDD}, that the entropy numbers of this unit ball decay like $n^{-s/d}$ and so in view of the Carl type inequality of Theorem \ref{mapentthm}  we have
 \be
 \label{lbC}
 \delta_{n,\gamma}^*(K)_{L_p}\ge cn^{-s/d},\quad n\ge 1.
\ee
    The main question therefore is whether the inequality in \eref{lbC} can be reversed.  In the case $p=2$,  the fact that it can be reversed follows from Theorem \ref{TH2}. The situation for
 $p\neq 2$ is not so straightforward and is still not settled.  Let us remark that for the weaker notion of  manifold widths 
 $\delta_n(K)_{L_p}$ both \eref{lbC} and its reverse have been proven,   see Theorem 1.1 in \cite{DKLT}.
 \vskip .1in
 \noindent
 {\bf Acknowledgment:}   The authors thank Professor Giles Godefroy for insightful discussions on the results of this paper.

\vskip .1in
\noindent
{\bf Affiliations:}

\noindent
Albert Cohen,  Laboratoire Jacques-Louis Lions,
Sorbonne Universit\'e, 4, Place Jussieu, 75005 Paris, France,  cohen$@$ann.jussieu.fr

\vskip .1in
\noindent
Ronald A. DeVore, Department of Mathematics, Texas A$\&$M University, College Station, TX 77843, rdevore$@$math.tamu.edu.
\vskip .1in
\noindent
Guergana Petrova, Department of Mathematics, Texas A$\&$M University, College Station, TX 77843,  gpetrova$@$math.tamu.edu.
\vskip .1in
\noindent
Przemys{\l}aw Wojtaszczyk, Institut of Mathematics Polish Academy of Sciences, ul. 
{\'S}niadeckich 8,  00-656 Warszawa, Poland, wojtaszczyk$@$impan.pl
   \end{document}